\documentclass[a4paper]{amsart}

\usepackage{amssymb,latexsym}
\usepackage{graphics}

\theoremstyle{plain}
\newtheorem{lemma}{Lemma}[section]
\newtheorem{proposition}[lemma]{Proposition}
\newtheorem{theorem}[lemma]{Theorem}

\theoremstyle{remark}
\newtheorem{remark}[lemma]{Remark}

\theoremstyle{definition}
\newtheorem{definition}[lemma]{Definition}

\numberwithin{equation}{section}

\newcommand{\R}{\mathbb{R}}
\newcommand{\N}{\mathbb{N}}
\newcommand{\Z}{\mathbb{Z}}
\newcommand{\s}{\text{\rm span}}
\newcommand{\TV}{\text{\rm Tot.Var.}}

\begin{document}

\title[$L^\infty$ solutions of hyperbolic systems]{Stability of
$L^\infty$ Solutions for Hyperbolic Systems with Coinciding Shocks and Rarefactions}

\author{Stefano Bianchini}
\address{S.I.S.S.A. - I.S.A.S., via Beirut 2-4, 34014 TRIESTE (ITALY)}
\email{bianchin@sissa.it}
\urladdr{http://www.sissa.it/\textasciitilde bianchin/}
\keywords{Hyperbolic systems, conservation laws, well posedness}
\subjclass{35L65}
\date{June 10, 2000}
\thanks{I wish to thank Prof. Alberto Bressan for his careful reading
of the manuscript}

\begin{abstract}
We consider a hyperbolic system of conservation laws 
\[
\left\{ \begin{array}{c}
u_t + f(u)_x = 0 \\
u(0,\cdot) = u_0
\end{array} \right.
\]
where each characteristic field is either linearly degenerate or
genuinely nonlinear. Under the assumption of coinciding shock and
rarefaction curves and the existence of a set of Riemann coordinates $w$, we
prove that there exists a semigroup of solutions $u(t) = \mathcal{S}_t
u_0$, defined on initial data $u_0 \in L^\infty$. The semigroup
$\mathcal{S}$ is continuous w.r.t. time and the initial data $u_0$ in
the $L^1_{\text{loc}}$ topology. Moreover $\mathcal{S}$ is unique and
its trajectories are obtained as limits of wave front tracking approximations.
\end{abstract}

\maketitle

%\vskip 1cm
\centerline{S.I.S.S.A. Ref. 65/2000/M}
\vskip 1cm

\section{Introduction}

Consider the Cauchy problem for a strictly hyperbolic system of
conservation laws 
\begin{equation}\label{E:hcl0}
\left\{ \begin{array}{c}
u_t + f(u)_x = 0 \\
u(0,\cdot) = u_0
\end{array} \right.
\end{equation}
where $u \in \R^n$ and $f : \Omega \mapsto \R^n$ is sufficiently
smooth, $\Omega$ open.
If the initial data $u_0$ is of small total variation, the global
existence was proved first in \cite{Gl}. Moreover a series of paper
\cite{Br1,Br2,BrCP,BrCo,BrLY} establishes the uniqueness and well
posedness of the Cauchy problem \eqref{E:hcl0}. However, when $u_0$ has
large total variation or even more generally $u_0$ belongs to $L^\infty$, the
solution $u$ may not exist globally in $L^\infty$ \cite{Jen}: 
only for special system it is possible
to consider initial data with large total variation. We recall some of
the results available in this direction.
\begin{itemize}
\item[1)] For scalar conservation laws, the entropic solution to
\eqref{E:hcl0} generates a contractive semigroup w.r.t the $L^1$
distance, on a domain of $L^\infty$ data \cite{Kr}. 
\item[2)] For general Temple class system, in \cite{BaBr,Bia,Ser1} it
is proved the existence and stability of the entropic solution for
initial data with arbitrarily large but bounded total variation.
\item[3)] If all characteristic families are genuinely nonlinear
and the system is Temple class, the existence and stability for
initial data in $L^\infty$ is proved
in \cite{BrGo}.
\item[4)] For special $2 \times 2$ systems, in which one of the
equation is autonomous, various results have been proved in
\cite{BaJe,BrSh}, with initial data with unbounded total variation.
\end{itemize}
An open question is if the semigroup of solutions to the systems of
case 2), defined on all the initial data $u_0$ with total variation
arbitrary large but bounded, can be extended to data in $L^\infty$. 
In many model systems, in fact, some of the characteristic fields are linearly
degenerate, so that the result in \cite{BrGo} does not apply: for
example, the traffic model considered in \cite{AwRa} or a $2 \times 2$
model for chromatography. Aim
of this paper is to prove that, at least in the case where the flux
function $f$ is convex and shocks and rarefactions coincide, the
solution to \eqref{E:hcl0} can be defined for $u_0 \in L^\infty$. 

To illustrate the heart of the matter, we assume that the system
\eqref{E:hcl0} admits a system of Riemann coordinates $w \in \R^n$,
and that shock and rarefaction curves coincide in $\Omega$. Moreover
we assume that each characteristic field is linearly degenerate or
genuinely nonlinear. We do not assume that rarefaction curves are
straight lines. We consider a set $E$ of the form
\[
E \doteq \Bigl\{ u \in \Omega: \ w(u) \in [a_i,b_i], \ i=1,\dots,n
\Bigr\}.
\]
With $L^\infty(\R;E)$ we denote the space of $L^\infty$ functions with
values in $E$. The main result of this paper is the following:
\begin{theorem}\label{T:final0}
There exists a unique semigroup $\mathcal{S} : [0,+\infty) \times L^\infty(\R;E)
\longmapsto L^\infty(\R;E)$ such that the following
properties are satisfied:
\begin{itemize}
\item[i)] for all $u_n,u \in L^\infty(\R;E)$, $t_n,t \in [0,+\infty)$,
with $u_n \to u$ in $L^1_{\text{loc}}$, $|t - t_n| \to 0$ as $n \to +\infty$, %[E:contin1]
\[
\lim_{n \to +\infty} \mathcal{S}_{t_n} u_n = \mathcal{S}_{t} u \quad
\text{in} \ L^1_{\text{loc}}; 
\]
\item[ii)] each trajectory $\mathcal{S}_t u_0$ is a weak entropic
solution to the Cauchy problem \eqref{E:hcl0}
with $u_0 \in L^{\infty}(\R;E)$;
\item[iii)] if $u_0$ is piecewise constant, then, for $t$ sufficiently
small, $\mathcal{S}_t u_0$ coincides with the function obtained by
piecing together the solutions of the corresponding Riemann problems.
\end{itemize}
\end{theorem}
The first two properties imply the well posedness of the Cauchy problem
\eqref{E:hcl0} with initial data in $L^\infty$. The uniqueness follows
because $\mathcal{S}$ satisfies iii) and it is limit of wave front
approximations. 

As it is shown in the
last example of \cite{BrGo}, the semigroup $\mathcal{S}$ cannot be uniformly
continuous: thus we cannot apply any compactness argument to construct
the solution $u(t) \doteq \mathcal{S}_t u_0$. The fundamental problem is
that, differently from \cite{BrGo}, the total variation of the Riemann
invariants corresponding to linearly degenerate families does not
decrease in time. 

The main idea of this paper is to study how the
solution to the characteristic equation 
\begin{equation}\label{E:charaeq0}
\dot x(t) = \lambda_i(u(t,x(t))), \quad x(0)=y,
\end{equation}
depends on the solution $u$ of \eqref{E:hcl0}. 
It will be shown that, for a fixed time
$\tau$, the map $y \mapsto x(\tau)$ depends Lipschitz continuously on
the initial data $u_0$, and moreover the Lipschitz constant is
independent of the total variation of $u_0$. Since the Riemann
invariant $w_i$ is the broad solution to
\[
(w_i)_t + \lambda_i(u(t,x)) (w_i)_x = 0,
\]
a simple argument gives the convergence of the wave front tracking
approximations. Note that this also implies the well posedness of the ODE
\eqref{E:charaeq0} when $u(t,x)$ is an $L^\infty$ solution of the
system \eqref{E:hcl0}. As noted in \cite{BrSh}, for general system
the solution to \eqref{E:charaeq0} does not exist. 

The paper is organized as follows.
Section \ref{S:setting} contains the basic assumptions on the system
\eqref{E:hcl0}. Moreover we construct the wave front approximation of
the solution $u(t)$. In Section \ref{S:shiftdif} we analyze carefully the
shift differential map, i.e. the evolution of a perturbation in
$u_0$ in which only the position of the initial jumps has
changed. The method we use is essentially the one in \cite{BrGo}, with
slight modifications due to the fact that in our system the rarefaction
curves do not need to be straight lines. The main result is here the
explicit computation of the shift differential map.

Section \ref{S:charact} is concerned with the equation for
characteristics \eqref{E:charaeq0}. We prove the Lipschitz dependence
of the map $y \mapsto x(t)$ w.r.t. both the initial data $u_0$ and
$y$. Moreover we will show that the Lipschitz constant is independent
from the total variation of $u_0$. Finally, in Section
\ref{S:semicont}, we prove Theorem \ref{T:final0}.

\section{Basic assumptions and wave front approximations}\label{S:setting}
%[S:setting]

We consider a strictly hyperbolic system of conservation laws %[E:hcl1]
\begin{equation}\label{E:hcl1}
u_t + f(u)_x = 0,
\end{equation}
where $f:\Omega \to \R^n$ is a smooth vector field defined on some
open set $\Omega \subseteq \R^n$. Let $A(u) \doteq DF(u)$ be the
Jacobian matrix of $f$ and denote with $\lambda_i(u)$ its eigenvalues and
with $r_i(u)$, $l^i(u)$ its right and left eigenvectors, respectively.
We assume that 
the eigenvalues $\lambda_i$ can be either genuinely nonlinear or
linearly degenerate. 
In the following the $i$-th rarefaction curve
through $u \in \Omega$ 
%i.e. the solution to the ODE
%\[
%\dot u = r_i(u), \qquad u(0)=u,
%\]
will be written as $R_i(s)u$, with $R_i(0)u=u$, while the $i$-th shock
curve will be denoted by $S_i(s)u$, and its speed by $\sigma_i(s,u)$. 
The directional derivative of a
function $\phi(u)$ in the direction of $r_i(u)$ will be denoted as
\[
r_i \bullet \phi(u) \doteq \lim_{h \to 0} \frac{\phi(u + h r_i(u)) -
\phi(u)}{h},
\]
while the left and right limit of a BV function $f$ in a point $x$
will be written as
\[
f(x-) = \lim_{y \to x-} f(y), \qquad f(x+) = \lim_{y \to x+} f(y).
\]
We assume that 
the rarefaction curves $R_i$ generate a system of
Riemann coordinates $w(u)$. 
We recall that a necessary and sufficient
condition for the local existence of Riemann coordinates is the
Frobenius involutive condition: if $[X,Y]$ denotes the Lie bracket of the
vector fields $X$,$Y$, the condition is
\[
[r_i,r_j] \in \s \{r_i,r_j\} \quad \text{for all}\ i,j =1,\dots,n.
\]
In the following we will use indifferently
the conserved coordinates $u$ or the Riemann coordinates $w$. 

Fix a domain %[E:domain]
\begin{equation}\label{E:domain}
E \doteq \Bigl\{ u \in \Omega: \ w(u) \in [a_i,b_i], \ i=1,\dots,n
\Bigr\}.
\end{equation}
Since $E$ is compact, there is a constant $c>0$ such that %[E:unifgnl]
\begin{equation}\label{E:unifgnl}
r_i \bullet \lambda_i(u) > c \quad \forall u \in E, \quad \text{if
}\lambda_i\text{ is genuinly nonlinear}.
\end{equation}
We suppose that the system \eqref{E:hcl1} is uniformly
strictly hyperbolic in $\Omega$: this means that there exists a
constant $d$ such that %[E:stricthyp]
\begin{equation}\label{E:stricthyp}
\lambda_{i+1}(u) - \lambda_i(v) \geq d, \quad \forall u,v \in E, \
i=1,\dots,n-1.
\end{equation}
We also assume that in the system \eqref{E:hcl1}
shock and rarefaction curves coincide: this implies \cite{Tem}
that either the rarefaction curve $R_i(s)u$ is a straight line or
the eigenvalue is linearly degenerate. 
In fact, one can prove that
\begin{equation}\label{E:eigendep}
\frac{d^2}{ds^2}\sigma_i(s,u)\Biggr|_{s=0} = \frac{1}{6} (r_i \bullet
\lambda_i(u)) \langle l^i(u), 
r_i \bullet r_i(u) \rangle + \frac{1}{3} r_i \bullet (r_i \bullet
\lambda_i(u)),
\end{equation}
and for the shock curve $S_i(s)u$ we have %[E:shockdep]
\begin{equation}\label{E:shockdep}
\langle l^j(u), S'''(0)u - R'''(0)u \rangle = \frac{1}{2(\lambda_j(u) -
\lambda_i(u))} (r_i \bullet \lambda_i(u)) \langle l^j(u), r_i \bullet
r_i(u) \rangle.
\end{equation}
If $\lambda_i$ is genuinely nonlinear, the left hand
side of \eqref{E:shockdep} is zero if and only if the
rarefaction curve is a straight line, 
because $r_i \bullet r_i(u)$ is orthogonal to $r_i(u)$.

The flux function $f$ thus satisfies the following assumptions:
\begin{itemize}
\item[{\bf H1)}] the eigenvalues $\lambda_i$ of $Df$ are linearly degenerate
or genuinely nonlinear;
\item[{\bf H2)}] the rarefaction curves form a system of coordinates;
\item[{\bf H3)}] shock and rarefaction curves coincide.
\end{itemize}
The system \eqref{E:hcl1} has thus $n_{ld}$ linearly degenerate
fields $\lambda_{i}$, corresponding to the Riemann invariants
$w_{i}$, and $n_{gnl} = n - n_{ld}$ genuinely nonlinear
fields $\lambda_{k}$, corresponding to the Riemann invariants
$w_{k}$. In the latter case we have $r_{k} \bullet
r_{k}(u) = 0$ for all $u \in E$. %[R:2x2case]

\begin{figure}
\centerline{\resizebox{16cm}{7cm}{{\includegraphics{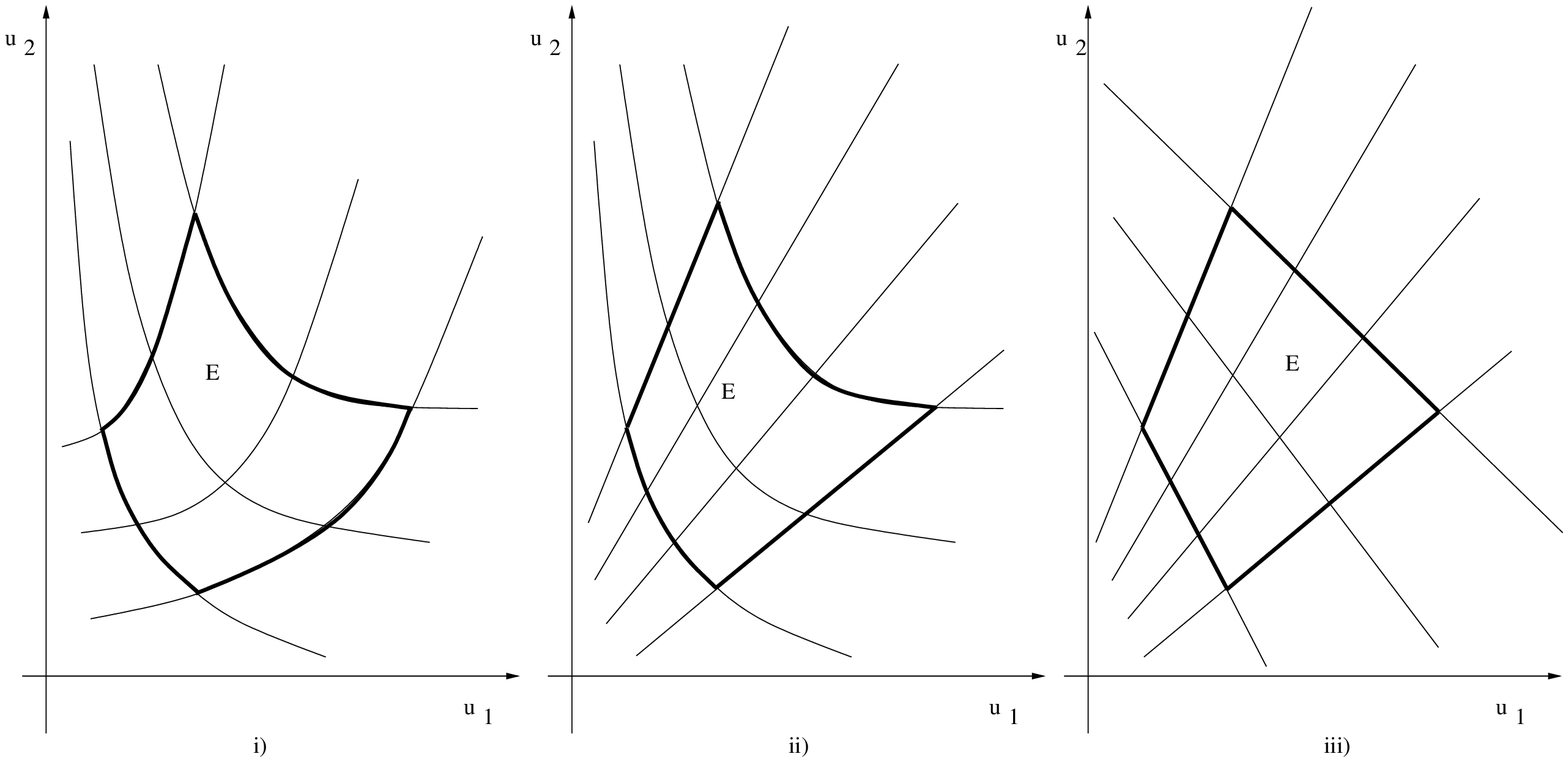}}}}
\caption{The various situations for a $2 \times 2$ system considered
in Remark \ref{R:2x2case}. %[Fi:2x2exmpl]
} \label{Fi:2x2exmpl}
\end{figure}

\begin{remark}\label{R:2x2case}
If $\Omega \subseteq \R^2$, then the rarefaction curves $R_i(s)u$
always generate a system of Riemann coordinates. Thus our assumptions
are satisfied by the following classes of systems:
\begin{itemize}
\item[i)] both eigenvalues are linearly degenerate;
\item[ii)] one eigenvalue is linearly degenerate, the other
genuinely nonlinear and the rarefaction curves of the latter are straight lines;
\item[iii)] both eigenvalues are genuinely nonlinear and the system is
of Temple class.
\end{itemize}
The various situations are shown in fig.~\ref{Fi:2x2exmpl}. Case ii)
corresponds to the traffic model considered in \cite{AwRa}, 
while case i) corresponds to $2 \times 2$ chromatography.
\end{remark}

Given the two points $u^-, u^+ \in E$, with coordinates
$u^- = u(w_1^-,\dots,w_n^-)$ and $u^+ = u(w_1^+,\dots,w_n^+)$, with
$w_i^+ \not = w_i^-$,
consider the intermediate states $u(\omega_i)$, where %[E:intersts]
\begin{equation}\label{E:intersts}
\omega_0= w(u^-), \quad \omega_i = (w_1^+,\dots,w_i^+,w_{i+1}^-,\dots,w_n^-), \quad
i=1,\dots,n.
\end{equation}
For all $i=1,\dots,n$, we denote with $v_i(u^-,u^+)$ the vectors
defined as %[E:linindep]
\begin{equation}\label{E:linindep}
v_i(u^-,u^+) = u(\omega_i) - u(\omega_{i-1}),
\end{equation}
and we define $r_i(u^-,u^+)$ as %[E:base1]
\begin{equation}\label{E:base1}
r_i(u^-,u^+) = 
\begin{cases} 
{\displaystyle \frac{v_i(u^-,u^+)}{|v_i(u^-,u^+)|}=\frac{u(\omega_i) -
u(\omega_{i-1})}{|u(\omega_i) - u(\omega_{i-1})|}} & \text{if } w_i^-
\not= w_i^+ \\
r_i(\omega_{i-1}) = r_i(\omega_i) & \text{if } w_i^- = w_i^+
\end{cases}
\end{equation}
where $r_i(u)$ is the $i$-th eigenvector of $DF(u)$. We assume that the vectors
$r_i(u^-,u^+)$ are linearly independent for all $u^-,u^+ \in E$. 
This condition is satisfied for data in a sufficiently
small neighborhood of a given point $\bar u \in \Omega$.
We denote also with
$\{ l^i(u^-,u^+), i=1,\dots,n\}$ the dual base.  
%Moreover let $P(u^-,u^+)$ be the projection
%operator on $\s \{r_i(u^-,u^+), i=1,\dots,n\}$.

We now define an approximated semigroup of solutions $\mathcal{S}^\nu$
on a set $E^\nu \subseteq E$. The construction is similar to the one in \cite{BaBr}.
For any integer $\nu \in \N$, set %[E:approxE]
\begin{equation}\label{E:approxE}
E^\nu \doteq \Bigl\{ u \in E: w_i(u) \in 2^{-\nu} \Z, i=1,\dots,n
\Bigr\}, 
\end{equation}
and let $D^{\nu,M}$ be the domain defined as %[E:approxdom]
\begin{equation}\label{E:approxdom}
D^{\nu,M} \doteq  \Bigl\{ u : \R \longmapsto E^{\nu}: u \ \text{piecewise
constant and } \TV(u) \leq M \Bigr\}. 
\end{equation}
Given $\bar u \in E^{\nu}$, we construct a solution $u(t)$ by wave
front tracking. We first define how to solve the Riemann problem
$[u^-,u^+]$, with $u^-, u^+ \in E^\nu$.

The solution to the Riemann problem $u^-,u^+$ is constructed by
piecing together the solutions to the simple Riemann problems
$[\omega_{i-1}, \omega_{i}]$, where $\omega_i$ is defined in
\eqref{E:intersts}. If the $i$-th field is linearly
degenerate, then $[\omega_{i-1}, \omega_{i}]$ is solved by a contact discontinuity
travelling with speed $\lambda_i(\omega_{i})$. If the $i$-th field
is genuinely nonlinear and $w_i^+ < w_i^-$, then $[\omega_{i-1},
\omega_{i}]$ is solved by a shock travelling with the Rankine-Hugoniot
speed $\sigma_i(\omega_{i-1},\omega_i)$. Finally, if the $i$-th field
is genuinely nonlinear and $w_i^+ > w_i^-$, then $[\omega_{i-1},
\omega_{i}]$ is solved by a rarefaction fan: if $w_i^+ = w_i^- + p_i
2^{-\nu}$, $p_i \in \N$, consider the states 
\[
\omega_{i,0} = \omega_{i-1}, \quad \omega_{i,l} = (w_1^+,\dots,w_{i-1}^+,w_i^- + \ell
2^{-\nu},w_{i+1}^-,\dots,w_n^-), \quad \ell=1,\dots,p_i.
\]
The solution will consist of $p_i$ shock waves
$[\omega_{i,l-1},\omega_{i,l}]$, travelling with the corresponding
shock speed $\sigma_i(\omega_{i,l-1},\omega_{i,l})$.

At time $t = 0$ we solve the initial Riemann problems of $\bar
u$. Note that the number of wave fronts is bounded by $2^\nu
\cdot \TV(\bar u)$. When two or more fronts interact, we solve
again the Riemann problem they generate, and so on. It is easy to show
that at each interaction at least one of the following alternatives
holds:
\begin{itemize}
\item[i)] the number of waves decreases at least by $1$;
\item[ii)] the total variation of the solution $u(t)$ decreases
by $2^{1-\nu}$,
\item[iii)] the interaction potential $Q(t)$, defined as %[E:intpoten]
\begin{equation}\label{E:intpoten}
Q(t) \doteq \sum_{\alpha, \beta \ \text{approaching}} |\sigma_\alpha|
|\sigma_\beta| \leq M^2, 
\end{equation}
decreases by $2^{-\nu}$. We recall that two waves $\sigma_\alpha$,
$\sigma_\beta$ of the families $k_\alpha$, $k_\beta$, located at
points $x_\alpha$, $X_\beta$, are considered as {\it approaching} if
$x_\alpha < x_\beta$ and $k_\alpha > k_\beta$.
\end{itemize}
This implies that there are at most a finite number of interactions,
so that we can construct our approximate solution for all $t \geq
0$. Note that $\mathcal{S}_t^{\nu} u = u(t)$ is a semigroup of
solutions, but not entropic due to the presence of rarefaction fronts.

%To extend the semigroup to $L^\infty$ data, we note that 
If the $i$-th family is linearly degenerate, the $i$-th Riemann
coordinate $w_{i}(t,\cdot)$ of the solution can be constructed by solving the
semilinear system %[E:semilin1]
\begin{equation}\label{E:semilin1}
\left\{ \begin{array}{cc}
(w_{i})_t + \lambda_{i}(u(t,x)) (w_{i})_x = 0 \\
w_{i}(0,x) = w_{i}(x) 
\end{array} \right.
\end{equation}
Since $u$ is a piecewise constant solution, with a finite number of
jumps, the broad solution to \eqref{E:semilin1} is well defined \cite{Br3}: if we
denote with $x(t,y)$ the solution to the ODE %[E:semilin2]
\begin{equation}\label{E:semilin2}
\dot x = \lambda_{i}(u(t,x)), \qquad x(0) = y,
\end{equation}
then the solution to \eqref{E:semilin1} is given by %[E:solut1]
\begin{equation}\label{E:solut1}
w_{i}(t,x(t,y)) = w_{i}(y).
\end{equation}
In the following sections we will consider the dependence on the initial
data $u_0$ of the
genuinely nonlinear Riemann coordinates $w_{k}(t,\cdot)$ and the
map $h_{i}^t(y)$ defined as %[E:lindegmap0]
\begin{equation}\label{E:lindegmap0}
h_{i}^t(y) \doteq x_{i}(t,y),
\end{equation}
where $x_{i}(t,y)$ is the solution to \eqref{E:semilin1}.

\section{Estimates on the shift differential map}\label{S:shiftdif}
%[S:shiftdif]

In this section we prove some properties of the shift differential
map. These properties are closely related to the structure of
\eqref{E:hcl1}, i.e. the conservation form, 
the coinciding shock and rarefaction assumption,
which prevents the creation of shock when two jumps of the same family
collide, and the existence of Riemann invariants, which prevents the
creation of shock when two jumps of different families interact.

Consider a wave front solution $u(t,\cdot)$ of \eqref{E:hcl1}, and
assume that the initial datum 
$u(0,\cdot)$ has a finite number $N$ of jumps $\sigma_\alpha$,
located in $y_\alpha$: 
\[
u(0,x) = \sum_{\alpha=1}^N \sigma_\alpha \chi_{[y_\alpha,+\infty)}(x).
\]
If $\xi_\alpha$ is the shift
rate of the jump $\sigma_\alpha$, define $u^\theta(t,\cdot)$ as the
front tracking solution with initial datum %[E:initshif]
\begin{equation}\label{E:initshif}
u^\theta(0,x) = \sum_{\alpha=1}^N \sigma_\alpha
\chi_{[y_\alpha + \theta \xi_\alpha,+\infty)}(x).
\end{equation}
In the following, we will use the integral shift function, defined by
%[E:intgrshif] 
\begin{equation}\label{E:intgrshif}
v(t,x) \doteq \lim_{\theta \to 0} \left\{ - \frac{1}{\theta}
\int_{-\infty}^x u^\theta(t,y) - u(t,y) dy \right\}.
\end{equation}
If $u(t,\cdot)$ has a shock $\sigma_\beta$, located in $y_\beta$, and if
$\xi_\beta$ is its shift rate, it is 
clear that the following relation holds: %[E:shifshok]
\begin{equation}\label{E:shifshok}
\sigma_\beta \xi_\beta = v(t,y_\beta+) - v(t,y_\beta-).
\end{equation}

We first recall the following result in \cite{BrGo}, obtained using
the conservation equation: %[L:cons]: 
\begin{lemma}\label{L:cons}
Consider a bounded, open region $\Gamma$ in the $t$-$x$ plane. Call
$\sigma_{\alpha}$, $\alpha=1,\dots,N$, the fronts entering $\Gamma$
and let $\xi_{\alpha}$ be their shifts. Assume that the fronts
leaving $\Gamma$, say $\sigma_\beta'$, $\beta=1,\dots,N'$, are linearly
independent. Then their shifts $\xi_\beta'$ are uniquely determined
by the linear relation %[E:cons]
\begin{equation}\label{E:cons}
\sum_{\beta=1}^{N'} \xi_\beta' \sigma_\beta' = \sum_{\alpha=1}^N
\xi_\alpha \sigma_\alpha. 
\end{equation}
\end{lemma}

\begin{remark}
As observed in \cite{BrGo}, formula \eqref{E:cons} implies that the shift
rates of the outgoing fronts depend only on the shift rates of the
incoming ones, and not on the order in which these wave-fronts
interact inside $\Gamma$. In particular we can perform the following
operations, without changing the shift rates of the outgoing fronts:
\begin{itemize}
\item[{\bf O1)}] switch the order of which three or more fronts interact;
\item[{\bf O2)}] invert the order of two fronts at time $0$, if they
have zero shift rate. 
\end{itemize}
\end{remark}

The second lemma is concerned with a configuration where a sequence of
contact discontinuities interacts with a wave of another family. %[L:lindeg]

\begin{lemma}\label{L:lindeg}
Consider a family of parallel contact discontinuities $\sigma_\alpha$, $\alpha=1,\dots,N$
of the $i$-th linearly degenerate family and a single wave-front $\sigma$ of
the $k$-th family, $k \not= i$. Let $\xi_\alpha$ and $\xi$ be
their initial shifts, respectively, and let $\xi_\alpha'$, $\xi'$ be their shifts
after interaction. Assume that $\xi_\alpha = \bar \xi$ for
all $\alpha$. Then after the interactions all the shift rates $\xi_\alpha'$ of the $i$-th 
family have the same value $\bar \xi'$ and %[E:lindeg]
\begin{equation}\label{E:lindeg}
\xi_\alpha' = \bar \xi' = {\bar \xi (\bar \Lambda' - \Lambda) - \xi
(\bar \Lambda' - \bar \Lambda) \over \Lambda - \bar \Lambda}, \qquad
\xi' =  {\bar \xi (\Lambda' - \Lambda) - \xi
(\bar \Lambda' - \bar \Lambda) \over \Lambda - \bar \Lambda}, 
\end{equation}
where $\bar \Lambda$, $\Lambda$ and $\bar \Lambda'$, $\Lambda'$
are the speeds of the shocks $\sigma_\alpha$ and $\sigma$,
before and after interaction, respectively.
\end{lemma}

\begin{figure}
\centerline{\resizebox{14cm}{5cm}{{\includegraphics{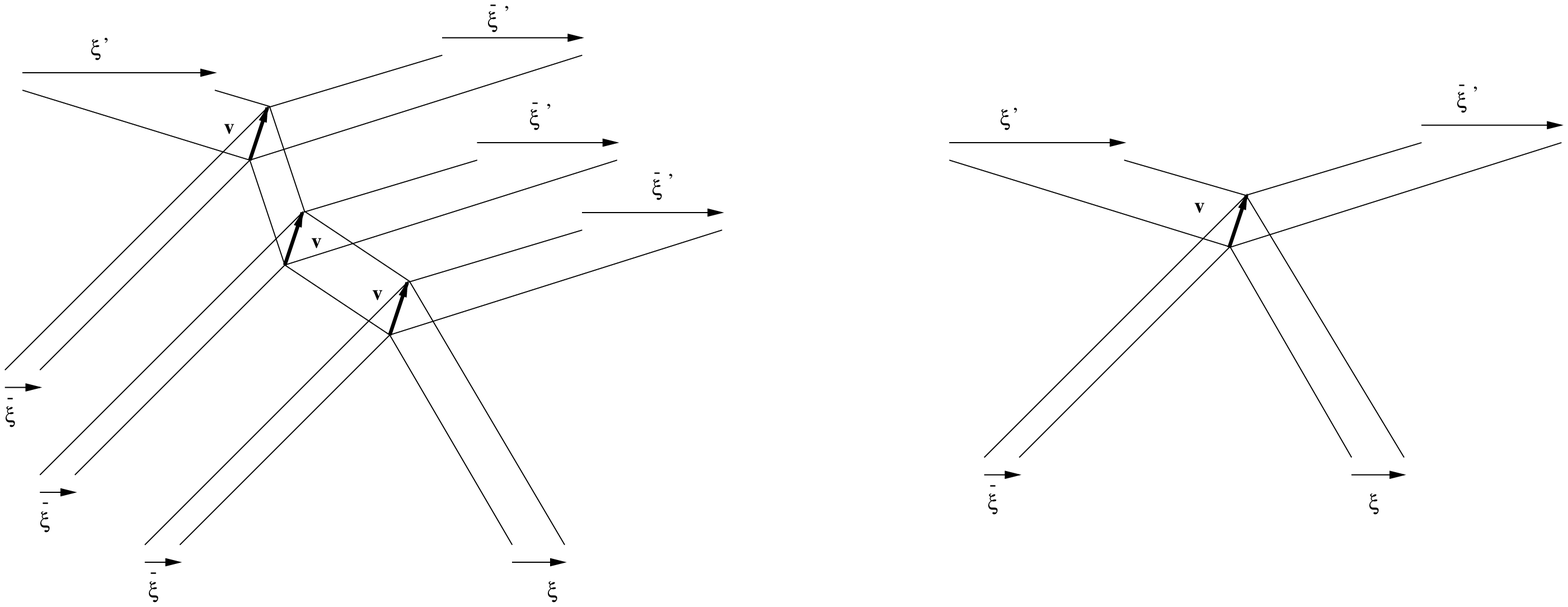}}}}
\caption{Interaction with a sheaf of shocks. %[Fi:Intlind]
} \label{Fi:Intlind}
\end{figure}

\begin{proof}
Define the vector $\mathbf{v}$ in the $t$-$x$ plane as the shift of
the first collision point. By a direct computation one finds %[E:vectv]
\begin{equation}\label{E:vectv}
\mathbf{v} = \left( {\xi - \bar \xi \over \bar \Lambda - \Lambda},
{\bar \Lambda \xi - \Lambda \bar \xi \over \bar \Lambda - \Lambda}
\right). 
\end{equation}
Since all the incoming shock of the linearly degenerate family have the
same speed $\bar \Lambda$, by simple geometrical considerations it
follows that the vector $\mathbf{v}$ is constant during all interactions
(fig.~\ref{Fi:Intlind}). Formula \eqref{E:lindeg} follows easily.
\end{proof}
\begin{remark}
Note that this lemma allows us to perform the following new operation,
without changing the shift rates:
\begin{itemize}
\item[{\bf O3)}] replace a family of contact discontinuities
$\sigma_\alpha$ of a linearly degenerate, all with the same shift rate
$\bar \xi$, by a single wave $\sigma = \sum \sigma_\alpha$ with shift rate $\bar \xi$.
\end{itemize}
\end{remark}

In the next lemma we will show that the existence of Riemann
coordinates $w$ implies a strong relation among shocks of different
families. %[L:Temple]

\begin{figure}
\centerline{\resizebox{14cm}{6cm}{{\includegraphics{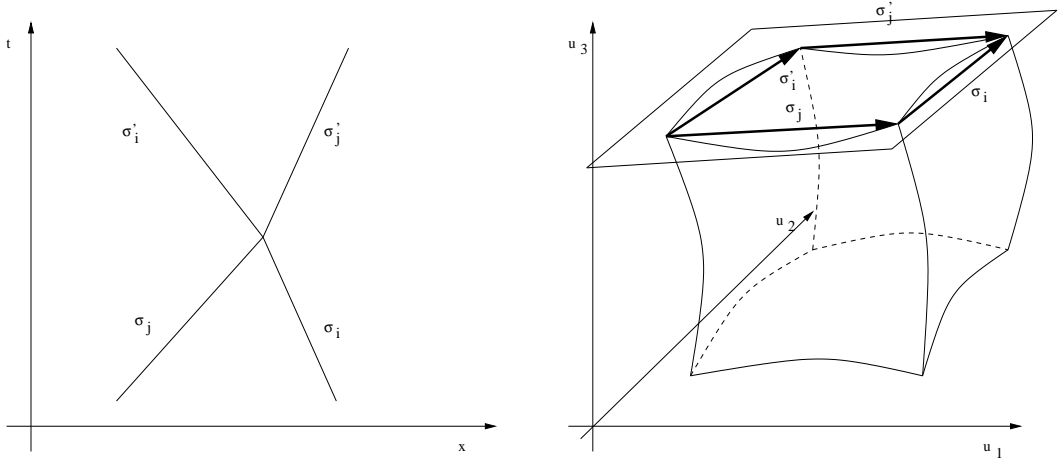}}}}
\caption{Vector relations among shocks. %[Fi:bidime]
} \label{Fi:bidime}
\end{figure}

\begin{lemma}\label{L:Temple}
Consider two adjacent jumps belonging to different families,
$\sigma_i$ and $\sigma_j$, $i < j$, located at $x_i > x_j$. Let
$\sigma_i'$, $\sigma_j'$ be their strength after interaction. Then
the following holds: %[E:Temple]
\begin{equation}\label{E:Temple}
\s \{ \sigma_i, \sigma_j\} = \s \{ \sigma_i', \sigma_j'\}.
\end{equation}
\end{lemma}
\begin{proof}
If $\xi_i$, $\xi_j$ are the shift rates before interaction, and
$\xi_i'$, $\xi_j'$ after interaction, then \eqref{E:Temple} follows
easily from the conservation relation %[E:consrel1]
\begin{equation}\label{E:consrel1}
\sigma_i \xi_i + \sigma_j \xi_j = \sigma_i' \xi_i' + \sigma_j'
\xi_j' \qquad \forall \xi_i, \xi_j \in \R, 
\end{equation}
because by assumptions no waves of other families are generated. By
condition \eqref{E:linindep} the conclusion follows. 
\end{proof}
%[R:O3]
\begin{remark}\label{R:O3}
Note that the previous lemma implies that the conservation relation
\eqref{E:consrel1} is bidimensional, i.e. the shocks $\sigma_i$,
$\sigma_j$ and $\sigma_i'$, $\sigma_j'$ lie on a two dimensional
plane (fig.~\ref{Fi:bidime}). 
We can obtain then an identity which relates the the strengths $\sigma$ with the speeds
$\Lambda$: substituting \eqref{E:lindeg} in
\eqref{E:consrel1}, since $\bar \xi$, $\xi$ are arbitrary, we get %[E:consrel2]
\begin{align}\label{E:consrel2}
\sigma_i (\Lambda_j - \Lambda_i) = \sigma_i' (\Lambda_i' -
\Lambda_j) + \sigma_j' (\Lambda_j' - \Lambda_j), \\
\sigma_j (\Lambda_i - \Lambda_j) = \sigma_i' (\Lambda_i' -
\Lambda_i) + \sigma_j' (\Lambda_j' - \Lambda_i). \notag
\end{align}
One can show that if a Riemann solver verifies \eqref{E:consrel2} for
all couple of waves $i,j$, then there
exists a flux function $f$ such that the wave front approximation is a
weak solution to \eqref{E:hcl1}.
%Note that \eqref{E:consrel2} is equivalent to the assumption that the
%system \eqref{E:hcl1} is in conservation form: if the Riemann solver
%satisfies \eqref{E:consrel2} for all couple of waves, then there
%exists a flux function $f$ such that he wave front approximation is a
%solution to \eqref{E:hcl1}.
\end{remark}

An important property of the shift differential map for Temple
class systems is the fact that a perturbation to the initial data, initially
localized in $[a,b]$, remains in the neighborhood of the set $\cup_i
[x_i(t,a),x_i(t,b)]$, where $x_i(t,y)$ is the solution of the $i$-th
characteristic equation starting at $y$. We now extend this property
to hyperbolic systems satisfying the hypotheses {\bf H1)}, {\bf H2)},
{\bf H3)} of section \ref{S:setting}. 

Consider $N$ jumps $\sigma_\alpha$, $\alpha=1,\dots,N$, of
some linearly degenerate family $i$, located at $x_\alpha$ and
corresponding to the jumps $c(\alpha) e_i$ in the Riemann coordinates $w$: %[E:invo1]
\begin{equation}\label{E:invo1}
\sigma_\alpha = u(w(x_\alpha-) + c(\alpha) e_{i}) - u(w(x_\alpha-)),
\end{equation}
for some constants $c(\alpha)$, $\alpha=1,\dots,N$.
%Consider $N$ shocks $\sigma_\alpha$, $\alpha=1,\dots,N$, of
%same linearly degenerate family $i_{ld}$, corresponding to the jumps
%$c(\alpha) e_i$ in the Riemann coordinates $w$: [E:invo1]
%\begin{equation}\label{E:invo1}
%\sigma_\alpha = u(w(\alpha-) + c(\alpha) e_{i_{ld}}) - u(w(\alpha-)),
%\end{equation}
%for some constants $c(\alpha)$, $\alpha=1,\dots,N$. 
%Consider the vectors [E:invol2]
%\begin{equation}\label{E:invol2}
%v_\alpha = u(w(\alpha)) - u(w(\alpha-1)), \qquad \alpha = 1,\dots,N,
%\end{equation}
%where $w(0) = w(1-)$.
\begin{definition}
We say that the jumps $\sigma_\alpha$ defined in \eqref{E:invo1} are
in involution if %[E:invol] 
\begin{equation}\label{E:invol}
\sum_{\alpha =1}^N c(\alpha) = 0,
\end{equation}
i.e. the initial and final Riemann coordinate $w_{i}$ is the
same: $w_{i}(x_1-) = w_{i}(x_N+)$. 
\end{definition} 
Note that, by the existence of Riemann coordinates, this relation does
not depend on the positions and strength of the shocks of the other
families. 
%holds for all $w(0) \in \{ w_{i_{ld}} = w_{i_{ld}}(0-) \} \cap E$, i.e. it
%does not depend on the positions of the shocks of the other families.
%\begin{figure}
%\centerline{\resizebox{12cm}{5cm}{{\includegraphics{lineardegen.eps}}}}
%\caption{Interaction with a sheaf of shocks. [Fi:Intlind]} \label{Fi:Intlind}
%\end{figure}
We can now extend Lemma 2 in \cite{BrGo} to our systems: %[L:excg]

\begin{lemma}\label{L:excg}
Consider a wave front tracking solution $u$.
Assume that there are $N$ shocks $\sigma_\alpha$
\begin{itemize}
\item[i)] either of the $i$-th linearly degenerate family in involution, 
\item[ii)] or of the $k$-th genuinely nonlinear family,
\end{itemize} 
and let $x_\alpha(t)$, $0 \leq t \leq T$, be the position of the shock
$\sigma_\alpha$, $\alpha=1,\dots,N$.
Then it is possible to assign shift rates to all shocks such that
$\xi_1 = 1$ and the shift of all fronts outside the strip $\Gamma
\doteq \{(t,x); t \in [0,T], x_1(t) \leq x \leq x_N(t) \}$ is zero. 
\end{lemma}

\begin{proof}
We consider only the case of linearly degenerate family $i$,
since in the other case the proof is exactly the one given in \cite{BrGo}. 

Let $x_\alpha(t)$, $\alpha = 1,\dots,N$, be the position of the
shock $\sigma_{\alpha}$ of the $i$-th family in involution, and
let $\bar w_{i}$ be the value of the Riemann coordinate at
$x_1(t)- = x_1(0)-$. For $w \in E$, define $\tilde w$ as
the projection of $w$ on the hyperplane $\{w_{i} = \bar
w_{i}\}$, and $\tilde u = u(\tilde w)$. 
%Call $\tilde
%r_{i_{ld}}(t,x)$ the unit vector in the direction of $u(t,x) -
%\tilde u(t,x)$: [E:direct1]
%\begin{equation}\label{E:direct1}
%\tilde r_i(t,x) \doteq \frac{u(t,x) - \tilde u(t,x)}{|u(t,x) - \tilde
%u(t,x)|}.
%\end{equation}
%If the latter is zero assume $\tilde r_i(t,x) = 0$.

We choose the shift rates such that %[E:choice]
\begin{equation}\label{E:choice}
- {d \over d\theta} \int_{-\infty}^x u^{\theta} dy = \sum_{x_i(t) \leq x}
  \xi(i)(t) \sigma_i(t) = c(t,x) (u(t,x) - \tilde u(t,x)),
\end{equation}
where $c(t,x)$ is a scalar function different from $0$ only in
$[x_1(t),x_N(t)]$, and we recall that $\tilde u(t,x) = u(\tilde w(t,x))$.
%where $\tilde r_i$ is the direction of the shock connecting the two
%states $w_i(0)$ and $w_i(t,x)$.

By imposing the value $\xi_1 =1$, i.e. $c(0,x_1(0)-) = 0$, $c(0,x_1(0)+) =
1$, we need to prove that
\eqref{E:choice} can be satisfied at time $t = 0$. We have two cases.
\begin{itemize}
\item[1)] If the jump $\sigma_i$ belongs to the $i$-th family and
is inside $[x_0(0),x_N(0)]$, then set $\xi = c(t,x_i-)$.
%\item[2)] If the jump $\sigma_i$ belongs to the $i_{ld}$-th family but
%it is not parallel to $\tilde r_i(t,x)$, then assign shift $c(0,x_\alpha-)$.
\item[2)] If the jump $\sigma_i$ belong to the $k$-th family with $k
\not= i$, then by assumption \eqref{E:linindep} and by \eqref{E:Temple}
there exists a unique shift $\xi_i$ and a unique constant $c(0,x+)$
such that
\[
\xi_i \sigma_i + c(0,x-)(u(0,x-) - \tilde u(0,x-)) = c(0,x+) (u(0,x+)
- \tilde u(0,x)). 
\]
\end{itemize}
Since we assume that the shocks are in involution, setting $\xi_N = c(0,x_n-)$
we have that \eqref{E:choice} holds at time $t=0$: in fact the last
jump has size $\tilde u(0,x_N(0)-) - u(0,x_N(0)-)$.

We now show that this property is conserved for all $t \geq
0$. This follows easily from conservation and Lemma
\ref{L:Temple}. The proof is exactly the same as in \cite{BrGo}: we
repeat it for completeness.
\begin{figure}
\centerline{\resizebox{12cm}{8cm}{{\includegraphics{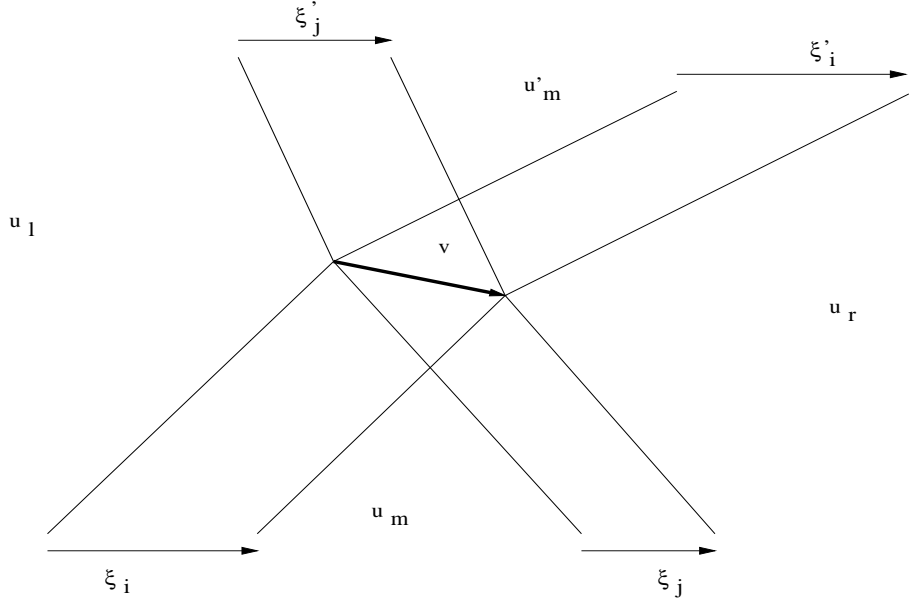}}}}
\caption{Computation of the shift rate. %[Fi:collision1]
} \label{Fi:collision1}
\end{figure}
Consider the interaction between two shocks $\sigma_i$ and $\sigma_j$
in the point $(\tau,y)$, see fig.~\ref{Fi:collision1}. By 
inductive assumption, we have for the states $u_l$, $u_m$ and $u_l$
that %[E:relatio1]
\begin{align}\label{E:relatio1}
\sum_{x_{\gamma}(\tau) < y} \sigma_{\gamma}(\tau) \xi_{\gamma}(\tau) &=
c_l (u_l - \tilde u_l), \\ 
c_l(u_l - \tilde u_l) + \sigma_i \xi_i &= c_m(u_m -
\tilde u_m), \notag \\ 
c_m(u_m - \tilde u_m) + \sigma_j \xi_j &= c_r(u_r - \tilde u_r). \notag
\end{align}
Using conservation we have %[E:relatio2]
\begin{equation}\label{E:relatio2}
\xi_i \sigma_i + \xi_j \sigma_j = \xi_j' \sigma_j' + \xi_i' \sigma_i',
\end{equation}
so that for the new middle state $u_m'$ we have %[E:relatio3]
\begin{equation}\label{E:relatio3}
c_l(u_l - \tilde u_l) + \sigma_j' \xi_j' = c_m' (u_m' - \tilde u_m') =
c_r(u_r - \tilde u_r) - \sigma_i' \xi_i', 
\end{equation}
and using Lemma \ref{L:Temple} we conclude
\[
\s \Bigl\{ u_l - \tilde u_l, \sigma_j'
\xi_j' \Bigr\} \bigcap \s \Bigl\{ u_r - \tilde u_r, \sigma_i'
\Bigr\} = \s \bigl\{u_m - \tilde u_m'\bigr\}.  
\]
The same relation proves that they vanish outside $\Gamma$: in fact,
assume for example that $c_l = 0$ and $j < i$. Then from
\eqref{E:relatio3} we get 
\[
\sigma_j' \xi_j' = c_m' (u_m' - \tilde u_m'), 
\]
which implies that $c_m' = 0$. This concludes the proof.
\end{proof}
%[R:lidegsh]
\begin{remark}\label{R:lidegsh}
%For waves of a linearly degenerate field $i_{ld}$,
%consider a wave $\sigma_{\alpha}$ of a linearly degenerate family. It
%is clear that we can think of this wave as the sum of wave of the same
%family $\sigma_\beta$, such that $\sum_\beta \sigma_\beta =
%\sigma_\alpha$. 
Note that for discontinuities of a linearly degenerate family all shift rates
has the same sign. 
Note moreover that if no waves of other families are present, then we shift
all jumps $\sigma_\alpha$ by unit rate $1$. This corresponds to the
case considered in Lemma \ref{L:lindeg}, i.e. to the substitution of a
family of contact discontinuities with a single jump, whose strength in this case is $0$
by the involution assumption.
\end{remark}

Using conservation and the previous lemmas, we obtain explicitly the
shift differential map at a given time $\tau$. 
%We consider first the case in which only one
%shock $\bar \sigma$ of the $k$-th family is shifted.
We recall that, given the states $u^-,u^+ \in E$, we denote with $r_i(u^-,u^+)$ the
vectors defined in \eqref{E:base1},  
and with $l^i(u^-,u^+)$ its dual base. 
Let $P_j(u^-,u^+)$ be the projection 
operator on $\s \{r_i(u^-,u^+), i=1,\dots,j\}$: %[E:projop]
\begin{equation}\label{E:projop}
P_j(u^-,u^+)v \doteq \sum_{i=1}^j \bigl\langle l^i(u^-,u^+), v
\bigr\rangle r_i(u^-,u^+),
\end{equation}
where $\langle \cdot, \cdot \rangle$ denotes the the scalar product in
$\R^n$. 

Given a point $(t,x)$, with $u(t,x)$ continuous in $x$,
define $x_i$ the intersection of the backward $i$-th characteristics
starting at $(t,x)$ with the real axis $\{(0,x)\}$, and for all
$(0,y)$ let $j(y)$ the index such that $x_{j(y)} \leq y < x_{j(y)-1}$,
$j(y) = 1,\dots,n+1$, with $x_{0} = +\infty$ and 
$x_{n+1} = -\infty$. Without any loss of generality, we can assume
that in $(0,y)$ there is a jump $\sigma$ of the $k$-th family.

Define the points $w_l,w_r \in E$ by %[E:riemnshft]
\begin{align}\label{E:riemnshft}
w_l(x,y) &\doteq 
\begin{cases}
w(t,x) & j(y) = 1\\
(w_1(0,y-),\dots,w_{j(y)-1}(0,y-),
w_{j(y)}(t,x),\dots,w_n(t,x)) & 2 \leq j(y) \leq n \\
w(0,y-) & j(y) = n+1
\end{cases} \\
%w_m(x,y) &\doteq 
%\begin{cases}
%w(0,y+) & j(y) = 1\\
%(w_1(t,x),\dots,w_k(0,y+),\dots,w_{j(y)-1}(t,x),
%w_{j(y)}(0,y),\dots,w_n(0,y)) & 2 \leq j(y) \leq n \\
%w(t,x) & j(y) = n+1
%\end{cases} \notag \\
w_r(x,y) &\doteq 
\begin{cases}
w(0,y+) & j(y) = 1\\
(w_1(t,x),\dots,w_{j(y)-1}(t,x),
w_{j(y)}(0,y+),\dots,w_n(0,y+)) & 2 \leq j(y) \leq n \\
w(t,x) & j(y) = n+1
\end{cases} \notag
\end{align}
Moreover define the point $w_m \in E$ by %[E:riemnshft2]
\begin{equation}\label{E:riemnshft2}
w_m(x,y) \doteq 
%\begin{cases}
%w(0,y+) & j(y) = 1\\
(w_1(t,x),\dots,w_k(0,y+),\dots,w_{j(y)-1}(t,x),
w_{j(y)}(0,y+),\dots,w_n(0,y+)) \quad 2 \leq j(y) \leq n, %\\
%w(t,x) & j(y) = n+1
%\end{cases} 
\end{equation}
if $k < j(y)$, and %[E:riemnshft1]
%\begin{align}\label{E:riemnshft1}
%w_l(x,y) &\doteq 
%\begin{cases}
%w(x-) & j(y) = 1\\
%(w_1(0,y-),\dots,w_{j(y)-1}(0,y-),
%w_{j(y)}(t,x),\dots,w_n(t,x)) & 2 \leq j(y) \leq n \\
%w(y-) & j(y) = n+1
%\end{cases} \\
%w_m(x,y) &\doteq 
%\begin{cases}
%w(y+) & j(y) = 1\\
%(w_1(t,x),\dots,w_{j(y)-1}(t,x),
%w_{j(y)}(0,y+),\dots,w_n(0,y+)) & 2 \leq j(y) \leq n \\
%w(x+) & j(y) = n+1
%\end{cases} \\
%w_r(x,y) &\doteq 
%\begin{cases}
%w(y+) & j(y) = 1\\
%(w_1(t,x),\dots,w_{j(y)-1}(t,x),
%w_{j(y)}(0,y+),\dots,w_n(0,y+)) & 2 \leq j(y) \leq n \\
%w(x+) & j(y) = n+1
%\end{cases}
%\end{align}
in a similar way, if $k \geq j(y)$, %[E:riemnshft3]
\begin{equation}\label{E:riemnshft3}
%w_m(x,y) &\doteq 
%\begin{cases}
%w(0,y+) & j(y) = 1\\
%(w_1(t,x),\dots,w_k(0,y+),\dots,w_{j(y)-1}(t,x),
%w_{j(y)}(0,y),\dots,w_n(0,y)) & 2 \leq j(y) \leq n \\
%w(t,x) & j(y) = n+1
%\end{cases} \\
w_m(x,y) \doteq 
%\begin{cases}
%w(0,y+) & j(y) = 1\\
(w_1(t,x),\dots,w_{j(y)-1}(t,x),
w_{j(y)}(0,y+),\dots,w_k(0,y-),\dots,w_n(0,y+)) \quad 2 \leq j(y) \leq n. 
%w(t,x) & j(y) = n+1
%\end{cases} 
\end{equation}
%If $r_i(w_1,w_2)$, $i=1,\dots,n$, are the vectors
%defined in \eqref{E:base1}, with $w_1,w_2 \in E$, let
%$P_i(w_1,w_2)$ be the projection operator on the subspace $\s
%\{r_k(w_1,w_2)\}$. 
%[E:base2]
%\begin{equation}\label{E:base2}
%r_i(x,y) = 
%\begin{cases}
%{\displaystyle \frac{v_i(w_l(x,y),w_r(x,y))}{|v_i(w_l(x,y),w_r(x,y))|}} & \text{if }
%w_{l,i}(x,y) \not= w_{r,i}(x,y), \\
%r_i(\omega_i) & \text{if } w_{l,i}(x,y) \not= w_{r,i}(x,y).
%\end{cases}
%\end{equation}
Define $P(x,y)$ as the vector %[E:projop1]
%be the projection of $u(w_r(x,y)) - u(w_l(x,y))$ on $\s
%\{r_k(x,y); w_{l,k}(x,y) \not= w_{r,k}(x,y) \}$: [E:projop1]
\begin{equation}\label{E:projop1}
P(x,y) \doteq 
\begin{cases}
0 & j(y) = 1 \\
P_{j(y)-1}(w_l,w_m) \sigma + P_{j(y)-1}(w_m,w_r)\bigl( \sigma -
P_{j(y)-1}(w_l,w_m)\sigma \bigr)
%\displaystyle \sum_{k=1 \atop w_{l,k}
%\not= w_{r,k} }^{j(y)-1} v_
%= \sum_{i=1}^{j(y)-1} r_i(w_l,w_r) \otimes l^i(w_l,w_r) 
& 2 \leq j(y) \leq n+1, k<j(y) \\
P_{j(y)-1}(w_l,w_m) \bigl( P_{j(y)-1}(w_m,w_r) \sigma \bigr) & 2 \leq
j(y) \leq n, k \geq j(y)
\end{cases}
\end{equation}
where $w_l = w_l(x,y)$, $w_m = w_m(x,y)$, $w_r = w_r(x,y)$ and 
$\sigma$ is the initial jump in $(0,y)$.
Consider now a front tracking solution $u^\theta$, obtained by
shifting the initial jumps $\sigma_\alpha$ in $y_\alpha$ with rates
$\xi_\alpha$. 
%We can then prove the following theorem: [T:shiftmap]
\begin{theorem}\label{T:shiftmap}
If $v(t,x)$ is the integral shift function of $u^\theta(t,\cdot)$,
defined in \eqref{E:intgrshif}, then %[E:shifttau]
\begin{equation}\label{E:shifttau}
v(t,x) = \lim_{\theta \to 0} \left\{ 
-\frac{1}{\theta}\int_{-\infty}^{x} u^\theta(t,y) - u(t,y) dy
\right\}  = \sum_{\alpha} P(x,y_\alpha) \xi_\alpha.
\end{equation}
\end{theorem}
\begin{proof}
The theorem will be proved outside the times of interaction, because
the Lipschitz dependence in $L^1$ of the approximate semigroup implies
the validity of \eqref{E:shifttau} for all $t \geq 0$.

If is sufficient to show that $\sum_{y_\alpha} P(x,y_\alpha)
\xi_\alpha$ is piecewise
constant, with jumps only at the points $x_\beta$ where $u(t,\cdot)$ has a
shock $\sigma_\beta$, and the following relation holds: %[E:rel1]
\begin{equation}\label{E:rel1}
\sum_{y_\alpha} \bigl( P(x_\beta+,y_\alpha) - P(x_\beta-,y_\alpha) \bigr)
\xi_\alpha = \sigma_\beta \xi_\beta, \qquad \lim_{x \to
-\infty} \sum_{y_\alpha} P(x,y_\alpha) \xi_\alpha = 0,
\end{equation}
where $\xi_\beta$ is the shift rate of $\sigma_\beta$, located in
$x$. Note that by \eqref{E:projop1} the second equality of
\eqref{E:rel1} is trivially satisfied.
\begin{figure}
\centerline{\resizebox{14cm}{5cm}{{\includegraphics{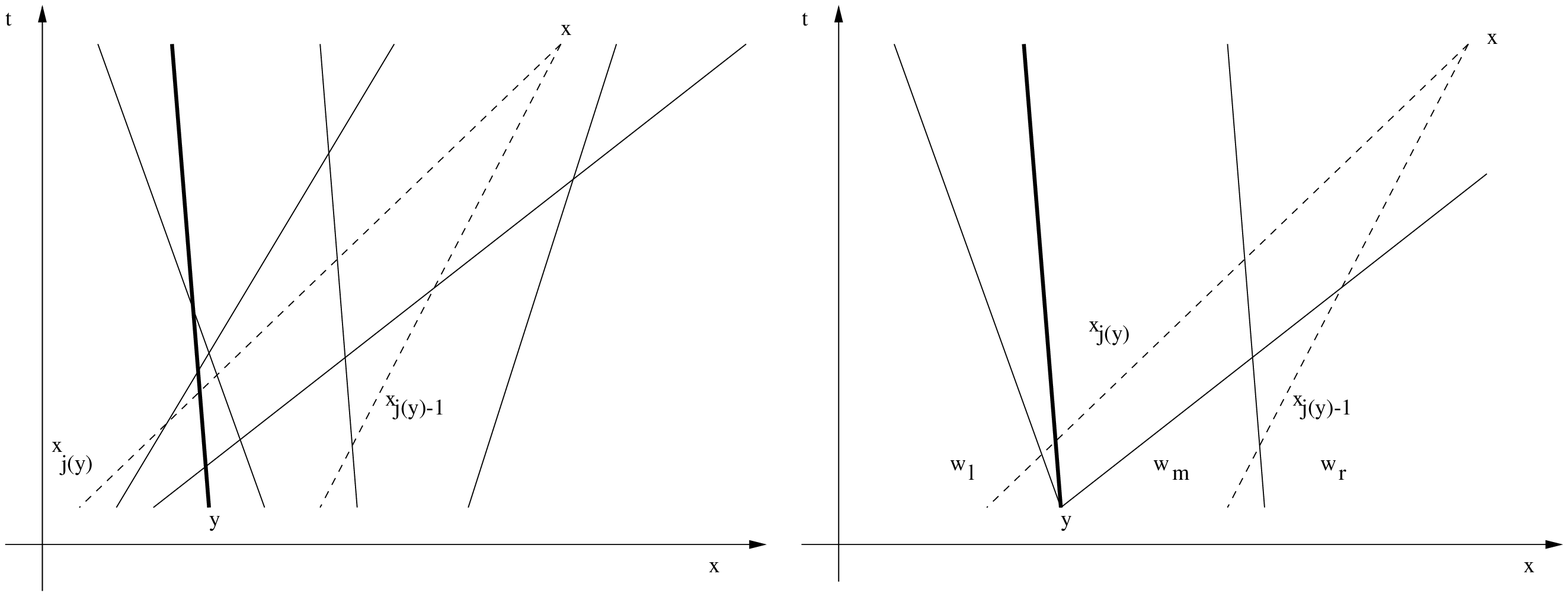}}}}
\caption{Wave pattern for the computation of formula \ref{E:projop}. %[Fi:riemnprb]
}\label{Fi:riemnprb}
\end{figure}

By linearity in the shift rates $\xi_\alpha$, we can consider the case
in which a single shock is shifted, let us say $\sigma$ at $y$: \eqref{E:shifttau}
becomes %[E:shiftsimple]
\begin{equation}\label{E:shiftsimple}
v(t,x) = P(x,y) \xi.
\end{equation}
Formula \eqref{E:projop} follows from the following considerations:
consider a wave front pattern, fig.~\ref{Fi:riemnprb}, where for
simplicity we assume that $k < j(y)$. 
The states $w_l$, $w_m$ are computed considering the Riemann problem
generated by adding to the $k$-jump $\sigma$ in $(0,y)$ all the
$i$-waves starting from the left of $(0,y)$ and ending in the right of
$(t,x)$ and  all the $i$-waves, with $i \not= k$, 
starting from the right of $(0,y)$ and ending in the left of
$(t,x)$. The jump $w_m,w_r$ is a single wave of the $k$-th
family formed by adding all the $k$-waves between $(0,y)$ and
$(t,x)$. Using the definition of $v(t,x)$ given \eqref{E:intgrshif}, one
obtains easily the second case of \eqref{E:projop1}: in fact the shift
rates of the shocks in the left of $(t,x)$ is given by the shift rates
of the jumps of the Riemann problem $w_l,w_m$ ending in the left of
$(t,x)$, $P_{j(y)-1}(w_l,w_m)\sigma$, plus the shift rate of the shock
$w_m,w_r$, $P_k(w_m,w_r)(\sigma - P_{j(y)-1}(w_l,w_m)\sigma)$. Since
only the $i$-waves with $i \geq j(y) > k$ are present in $\sigma -
P_{j(y)-1}(w_l,w_m)\sigma$, then
$P_k(w_m,w_r)(\sigma - P_{j(y)-1}(w_l,w_m)\sigma) =
P_{j(y)-1}(w_m,w_r)(\sigma - P_{j(y)-1}(w_l,w_m)\sigma)$. The other
cases can be computed in a similar way: in this case one solves the
Riemann problem $w_m,w_r$ in $(0,y)$, and consider the $k$-wave
$w_l,w_m$ starting in the left of $(0,y)$ and ending in the right of
$(t,x)$. 

From the above considerations it is clear that $P(x,y)$ is piecewise
constant, with jumps only when in $(t,x)$ there is a $i$-shock $\sigma'$:
in fact otherwise the wave front pattern used to compute $P(x,y)$ remains the same.
%In fact if $w(t,x-) =
%w(t,x+)$ and $j(y)$ increases of $1$, then there
%are no $j(y)$-waves in $(0,y)$ and then
%$w_{j(y)}(t,x)=w_{j(y)}(0,y-)=w_{j(y)}(0,y+)$. As a consequence 
%$w_{l}(x-,y)=w_{l}(x+,y)$, $w_{m}(x-,y) = w_m(x+,y)$ 
%and $w_{r}(x-,y)=w_{r}(x+,y)$, and thus: 
%\[
%_{j(y)}(w_l,w_m) \sigma = P_{j(y)-1}(w_l,w_m)\sigma, \qquad
%P_{j(y)}(w_m,w_r) \sigma = P_{j(y)-1}(w_m,w_r) \sigma.
%\] 
%Notice first that if $k$ is genuinely nonlinear, then from $r_k
%\bullet r_k = 0$ it follows easily that [E:projop2]
%\begin{equation}\label{E:projop2}
%P(x,y) = 
%\begin{cases}
%0 & j(y) = 1 \\
%P_{j(y)-1}(w_l,w_r) \sigma 
%& 2 \leq j(y) \leq n+1, k<j(y) \\
%\end{cases}
%\end{equation}
%The same is true if there are no shocks of the $k$-th family between
%$(0,y)$ and $(t,x)$.
Let $\{z_p:p=1,\dots,M\}$ be the set of the starting points of all shocks
arriving in $(t,x)$, and define %[E:startpoi]
\begin{equation}\label{E:startpoi}
z^- = \min_p z_p, \qquad z^+ = \max_p z_p.
\end{equation}
We consider two cases:
\begin{itemize}
\item[1)] the shocks arriving in $(t,x)$ start on both sides of
$(0,y)$: $z^- \leq y \leq z^+$.
In this case, $(P(x+,y) - P(x-,y))\xi$ is the shift rate of the
$i$-shock starting in the Riemann problem $w_l(x-,y),w_m(x+,y)$ if $i > k$
($w_m(x-,y),w_r(x+,y)$ if $i < k$) which collides with a $k$-shock $w_m(x+,y),w_r(x+,y)$
($w_l,w_m$ if $i < k$): in fact the only difference is that in
$w_m(x-,y),w_m(x+,y)$ there is a shock of the $i$-th family starting
in $(0,y)$, and $i$ is genuinely nonlinear. Finally, using $r_i \bullet r_i(u) = 0$
and Lemma \ref{L:Temple}, one can change position to the $i$-wave and
the remaining $k$-wave $w_m,w_r$, whose strength does not change. %(fig.~\ref{Fi:riemnprb1}).  

If $i=k$, there are no
$k$-shocks starting on the right (left) of $(0,y)$ and ending on the
right (left) of $(t,x)$, so that $(P(x+,y) - P(x-,y))\xi$
is the shift rate of the $i$-shock of the Riemann
problem $w_l(x-,y),w_r(x+,y)$.
\item[2)] the shocks of the $i$-th family arriving in $(t,x)$ start
either in $(-\infty,y)$ or $(y,+\infty)$: assume for definiteness that
$y < z^-$. In this case the difference $(P(x+,y) - P(x-,y))\xi$ is the
shift rate of the shock $\sigma'$ colliding with the shifted shocks
of the Riemann problem $w_l(x-,y),w_m(x-,y)$ in $(0,y)$, crossing the jump
$w_m(x-,y),w_r(x-,y)$, and finally overtaking $\sigma'$. In fact one
can use Lemma \ref{L:Temple} (and $r_i \bullet r_i(u) = 0$ if $i$ is
genuinely nonlinear) to obtain the wave pattern of fig.~\ref{Fi:riemnprb2}.
\end{itemize}
%Thus, to verify \eqref{E:rel1}, we need to prove the
%following. Consider a shock $\sigma'$ of the $i$-th family, located at
%$(t,x)$. Let $\{z_p:p=1,\dots,M\}$ be the set of the starting
%points of all shocks arriving in $(t,x)$, and define
%[E:startpoi]
%\begin{equation}\label{E:startpoi}
%z^- = \min_p z_p, \qquad z^+ = \max_p z_p.
%\end{equation}
%Assume that 
%the jump $\sigma$, located in $(0,y)$ is shifted with
%rate $\xi$: then 
%the shift rate $\xi'$ of $\sigma'$ is
%\begin{equation}\label{E:diffproj}
%P(x+,y) - P(x-,y) = 
%\begin{cases} 
%P_{k<j}(x-) - P_{k<j}(x+) & \text{if } x_j(x) \leq y < x_{j-1}(x), j < i \\
%r_i(w_l,w_r) \otimes l^i(w_l,w_r) & \text{if } x_i(x-) \leq y \leq x_i(x+) \\
%P_{k>j}(x-) - P_{k>j}(x+) & \text{if } x_j(x) \leq y < x_{j-1}(x), j > i
%\end{cases}
%\end{equation}
%Then the shift rate of $\sigma$ is then given by 
%[E:shiftgen]
%\begin{equation}\label{E:shiftgen}
%\xi' \sigma' = 
%\begin{cases}
%0 & \text{if }j(y)=1,n+1 \\
%{\displaystyle \sum_{k=j(y)}^n r_k(x,y) \bigl\langle l^k(x,y), 
%u_r(x,y) - u_r(x,y) \bigl\rangle
%\bar \xi} & \text{if } j(y) > i  \\
%{\displaystyle \sum_{k=i+1}^n r_k(x,y) \bigl\langle l^k(x,y), 
%u_r(x,y) - u_r(x,y) \bigl\rangle
%\bar \xi} & \text{if } j(y) = i, y > z^+  \\
%r_i(x,y) \Bigl\langle l^i(x,y), u_r(x,y) - u_r(x,y) \Bigr\rangle
%\xi & \text{if } j(y)=i, z^- \leq \bar y \leq z^+ \\
%{\displaystyle \sum_{k=1}^{j(y)-1} r_k(x_\alpha,\bar y) \bigl\langle
%l^k(x,y), u_r(x,y) - u_l(x,y) \bigl\rangle \bar \xi} &
%\text{if } j(\bar y) \leq i, \bar y < z_\alpha^-  
%\end{cases}
%\end{equation}
The various cases will be proved in the following lemmas. 
%[L:simplify1]
\begin{lemma}\label{L:simplify1}
Assume that $z^- \leq y \leq z^+$, i.e. case 1). 
%Consider the states $w(\bar x(0)-)$,  $w(\bar x(0)+)$ and $w(\bar
%x(\tau)-)$, $w(\bar x(\tau)+)$, and construct the states [E:Riemprb1]
%\begin{align}\label{E:Riemprb1}
%w_l &= (w_1(0-),\dots, w_{k-1}(0-),w_k(\tau-),\dots,w_n(\tau-)), \\
%w_r &= (w_1(\tau+),\dots,w_k(\tau+),w_{k+1}(0+),\dots,w_n(0+)). \notag
%\end{align}
If the shock $\sigma'$ is of the $i$-th family, then its shift $\xi'$ is %[E:shiftsame]
\begin{equation}\label{E:shiftsame}
\xi' \sigma' 
%= r_i(w_l,w_r) \Bigl\langle l^i(w_l,w_r),
%u_r(x,y) - u_l(x,y) \Bigr\rangle \xi 
= \bigl( P(x+,y) - P(x-,y) \bigr)
\xi.
\end{equation}
\end{lemma}
\begin{figure}
\centerline{\resizebox{14cm}{5cm}{{\includegraphics{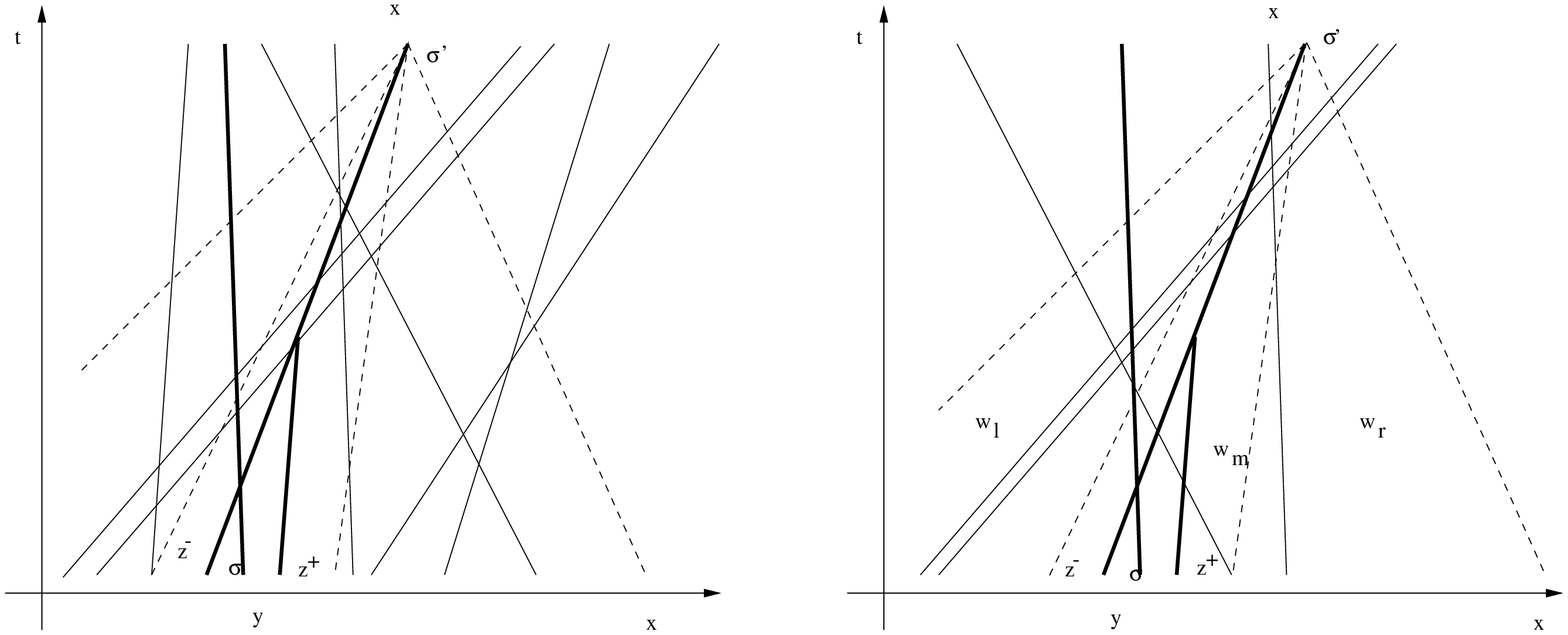}}}}
\caption{Computation of the shift rate in the case of Lemma
\ref{L:simplify1}. %[Fi:riemnprb1]
}\label{Fi:riemnprb1}
\end{figure} 
\begin{proof}
We follow closely the method of \cite{BrGo}. Assume for definiteness
$k < j(y)$, the other cases being similar.
The basic idea is to reduce the computation to the single Riemann
problem $w_l(x-,y),w_m(x+,y)$, with eventually a single $k$-wave
$w_m,w_r$. 

Consider fig.~\ref{Fi:riemnprb1}.
By Lemma \ref{L:cons}, we can simplify the wave
configuration considering only the fronts crossing
starting in the right of $(0,y)$ and ending in the 
left of $(t,x)$: in fact we can move the other fronts to $\pm \infty$
without changing the shift rate of $\sigma'$. 

We can now 
shift the initial position of the waves of the $i$-th family merging
in $x$ such that their initial position coincide with $y$,
without changing the shift rate $\xi'$. This operation can be
repeated for all shocks of genuinely nonlinear families.

Finally, we can move the shocks of the linearly degenerate
families such that they  have the same sequence of interaction with the other
shocks. This means that, if $x^j_{i}$ is the position of the
$j$-th shock of the $i$-th linearly degenerate family, the only interactions among
shocks occurring in the sector $[x^1_{i}(t), x_{i}^n(t)]$
are the one involving one $i$-th wave and one $k$-th wave, with
$k \not= i$. Using Lemma \ref{L:lindeg}, we can at this point
substitute them with a single shock, whose strength is the sum of the
strengths of the $i$-waves. Finally we move their position at $t
= 0$ such that it coincides with $y$: we obtain the wave patterns of
fig.~\ref{Fi:riemnprb1}. 
%\begin{figure}
%\centerline{\resizebox{12cm}{5cm}{{\includegraphics{simplify.eps}}}}
%\caption{[Fi:simplify]}\label{Fi:simplify}
%\end{figure} 
To conclude, we just need to prove that the Riemann problem obtained
in this way is exactly $w_l(x-,y),w_m(x+,y)$ and that the remaining
$k$-wave is $w_m(x+,y),w_r(x+,y)$.
%: in
%fact, if $\sigma_j$ is the strength of the $j$-th shock, by
%conservation it follows
%\[
%\xi \sigma = \xi' \sigma' +
%\sum_{j \not= i} \xi_{j} \sigma_{j}, 
%\]
%and \eqref{E:shiftsame} is a consequence of the linear independence of
%the vectors $v_i(w_l,w_r)$.

By the previous argument, the strength of the shock of the $j$-th
family $j < i$, $j \not= k$, is given by the $j$-waves starting in the right of
$y$ and ending in the left of $x$: since they are the only
$j$-wave crossing the segment $[(0,y-),(t,x+)]$, it follows
\[
w_{l,j}(x,y) = w_j(0,y-), \qquad
w_{r,j}(x,y) = w_j(t,x+).
\]
The other relations for $j = k$ and $j > i$ follows in the same way. 
Finally, for $j=i$ the jump is $w_i(t,x+) - w_i(t,x-)$. Note that the
wave pattern is the same obtained in 1).
\end{proof}
We consider only the case $y < z^-$, since the other is entirely similar.
%[L:simplify2]
\begin{lemma}\label{L:simplify2}
Assume that $y < z^-$. Then the shift $\xi'$ of $\sigma'$ is given by
%[E:shiftdiff1]
\begin{equation}\label{E:shiftdiff1}
\xi' \sigma' = 
%\sum_{k = j(y)}^{n} r_k(x-,y) \bigl\langle
%l^k(x-,y), u_r(x-,y) - u_l(x-,y) \bigr\rangle \\
%& \qquad \qquad - \sum_{k = j(y)}^{n} r_k(x+,y) \bigl\langle l^k(x+,y),
%u_r(x+,y) - u_l(x+,y) \bigr\rangle. \notag 
\bigl( P(x+,y) - P(x-,y) \bigr) \xi.
\end{equation}
\end{lemma}
\begin{proof}
The hypothesis implies that the $i$-th shocks ending at $x$
starts in the right of $y$. With the same simplification
considered in Lemma \ref{L:simplify1}, we reduce to the Riemann
problem $w_l(x,y),w_r(x,y)$ in $\bar y$, such
that the waves of the $j$-th families, $j > i$, generated at $\bar y$
collide with the $i$-wave in $x_\alpha$ (see fig.~\ref{Fi:riemnprb2}),
after overtaking the $k$-wave $w_m,w_r$. 
\begin{figure}
\centerline{\resizebox{14cm}{5cm}{{\includegraphics{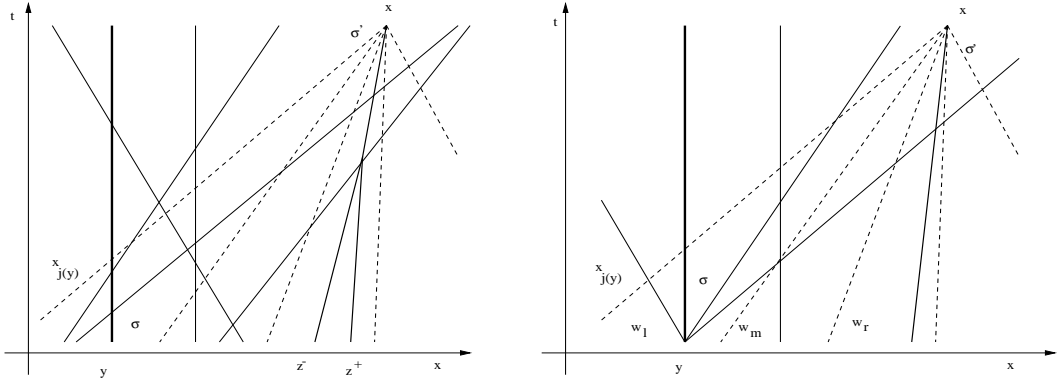}}}}
\caption{Computation of the shift rate in the case of Lemma
\ref{L:simplify2}. %[Fi:riemnprb2]
}\label{Fi:riemnprb2}
\end{figure} 
The conclusion follows easily, since the wave pattern is the same
considered in 2).
%The conclusion \eqref{E:shiftdiff1} follows by conservation: in fact,
%the amount of shift brought by the waves colliding with $\sigma'$ is given by
%\[
%\sum_{k = j(y)}^{n} r_k(x-,y) \bigl\langle
%l^k(x-,y), u_r(x-,y) - u_l(x-,y) \bigr\rangle \xi,
%\]
%while the outgoing waves have shift amount
%\[
%\xi' \sigma' + \sum_{k = j(y)}^{n} r_k(x+,y) \bigl\langle
%l^k(x+,y), u_r(x+,y) - u_l(x+,y) \bigr\rangle.
%\]
\end{proof}
This concludes the proof of Theorem \ref{T:shiftmap}.
\end{proof}
Finally we extend to our case the following result proved in \cite{BrGo}: %[P:finiteTV]
\begin{proposition}\label{P:finiteTV}
Let $u$ be a wave-front tracking solution, and consider two
wave-fronts, $x(t)$ and $y(t)$, $t \in [0,T]$. Then there exists a second front
tracking solution $\tilde u$ such that the initial and final positions of the two
shocks is the same, and $\TV(\tilde u)$ is uniformly bounded.
\end{proposition}

\begin{figure}
\centerline{\resizebox{15cm}{4cm}{{\includegraphics{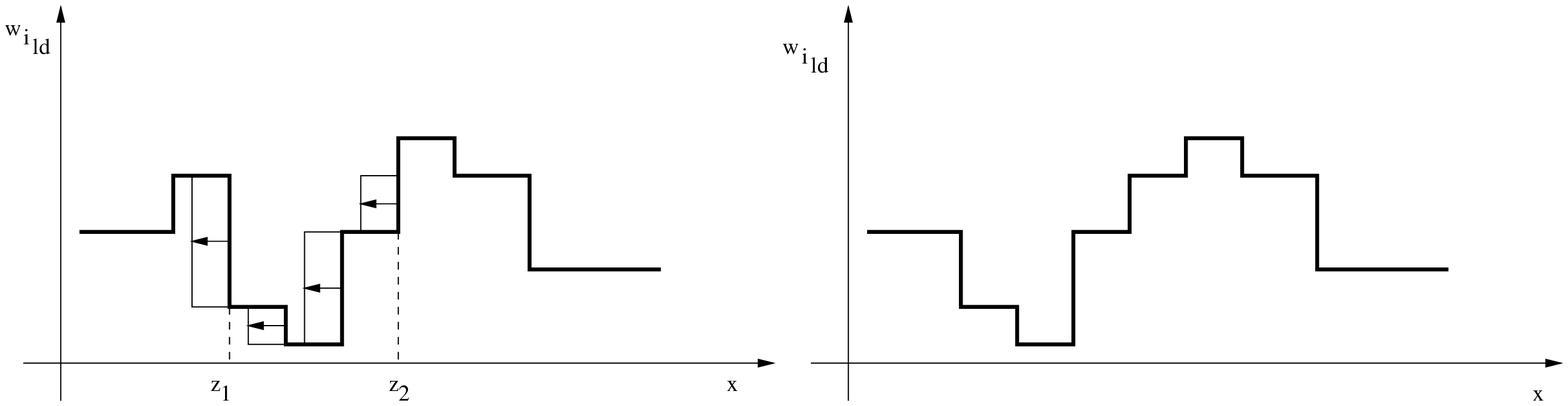}}}}
\caption{Cancellation among contact discontinuities. %[Fi:lidecanc]
}\label{Fi:lidecanc} 
\end{figure}

\begin{proof}
For genuinely nonlinear fields, the proof is the same as in \cite{BrGo}. We
then restrict the proof to the case of a linearly degenerate fields $i$.

Assume that there exists two jumps $\sigma_1$, $\sigma_2$ of the
$i$-th family, with positions $z_1(t) < z_2(t)$, such that %[E:regions1]
\begin{equation}\label{E:regions1}
x(0) \notin [z_1(0),z_2(0)] \ \text{and} \ y(0) \notin
[z_1(0),z_2(0)], \quad
x(T) \notin [z_1(T),z_2(T)] \ \text{and} \ y(T) \notin [z_1(T),z_2(T)].
\end{equation}
For definiteness, assume $w_{i}(0,z_1-) < w_{i}(0,z_1+)$, 
and the following conditions is satisfied:
%[E:cancrel1]
\begin{equation}\label{E:cancrel1}
w_{i}(0,z_1-) \in [w_{i}(0,z_2-), w(0,z_2+)].
%, \quad \text{or}
%\quad w_{i_{ld}}(z_2+) \in [w_{i_{ld}}(z_1-), w(z_1+)].
\end{equation}
Let $\sigma_\alpha$, $\alpha=1,\dots,N$ be the jumps of linearly degenerate family
$i$ in the strip $[z_1(0),z_2(0))$: if we define 
\[
\sigma_{N+1} = u(w_{i}(0,z_1-)) - u(w_{i}(0,z_2-)),
\]
it is easy to verify that the shocks $\sigma_\alpha$, $\alpha =
1,\dots,N+1$, are in involution. By Lemma \ref{L:excg}, we can then moving the jumps to the
left until either $z_1(t)$ meets the wave fronts $x(t)$, or $z_1(t)$
coincides with another shock of the $i$-th family (fig. \ref{Fi:lidecanc}). 
It is clear that we can repeat the same procedure also in the
following cases:
\begin{itemize}
\item[i)] $w_{i}(0,z_1-) > w_{i}(0,z_1+)$ and 
$w_{i}(0,z_1-) \in [w_{i}(0,z_2+), w(0,z_2-)]$;
\item[ii)] $w_{i}(0,z_2-) < w_{i}(0,z_2+)$ and 
$w_{i}(0,z_2+) \in [w_{i}(0,z_1+), w(0,z_1-)]$;
\item[iii)] $w_{i}(0,z_2-) > w_{i}(0,z_2+)$ and 
$w_{i}(0,z_2+) \in [w_{i}(0,z_1-), w(0,z_1+)]$.
\end{itemize}
It is now easy to prove that the total variation of the jumps of the
$i$-th family satisfying \eqref{E:regions1} can at most be 
$3 \|w\|_{\infty}$. Since $x(t)$, $y(t)$ divide the lines $t=0$ and
$t=\tau$ in three regions, the total variation of $w_{i}$ is
bounded by $27 \|w\|_{\infty}$.
\end{proof}

\section{Estimates on characteristics}
\label{S:charact}
%[S:charact]

In this section we prove some estimates on the solution $x_i(t,y)$ of the
characteristic equation: %[E:charaeq1]
\begin{equation}\label{E:charaeq1}
\left\{ \begin{array}{l}
\dot x_i = \lambda_{i}(u(t,x_i)) \\
x_i(0) = y
\end{array} \right.
\end{equation}
We assume for simplicity that the $i$-th family is linearly
degenerate, however the same results are valid for characteristics of a
genuinely nonlinear family if the following condition holds:
for all $\tau$ there exists an $\epsilon$ such that in
the strip $\{ (t,x); \tau \leq t \leq T, x_i(t,y)-\epsilon \leq x \leq
x_i(t,y)+\epsilon \}$ there are no shock waves of the $i$-th family. 
Given front tracking approximation $u$, $x_i(t,y)$ is
unique, since it crosses only a finite number of transversal
jumps, and it depends Lipschitz continuously on the initial data $y$
(see \cite{Br3}). 

We want to give uniform estimates on this dependence. The idea is to
suppose that in $y$ there is a shock $\sigma^\epsilon$ of the $i$-family of size
$\epsilon$: $w_{i}(0,y+) - w_i(0,y-) = \epsilon$. Since by
assumption no shocks of the $i$-family collide 
with $\sigma^\epsilon$, it is easy to construct a wave front solution:
for $x < x(t,y)$, the solution $u^\epsilon(t,\cdot)$ takes values in 
\[
E^{\nu,-} \doteq \Bigl\{ u : w_j(u) \in [a_j,b_j] \cap 2^{-\nu} \Z, j=1,\dots,n \Bigr\},
\]
while for $x > x(t,y)$, enlarging $E$ and assuming $\epsilon$
sufficiently small, 
\[
E^{\nu,+} \doteq \Bigl\{ u : w_j(u) \in [a_j,b_j] \cap 2^{-\nu} \Z, j
\not= i, w_i(u) \in [a_i,b_i] \cap \bigl\{ 
2^{-\nu} \Z + \epsilon \bigr\} \Bigr\}. 
\]

The following lemma proves the continuous dependence of the solution
$u^\epsilon(t)$ and the position $x^\epsilon_i(t,y)$ of the shock
$\sigma^\epsilon$ w.r.t. $\epsilon$.
%[L:littleshock]
\begin{lemma}\label{L:littleshock}
Consider a front tracking solution $u$, with initial data $u_0$ and the
characteristic lines $x_i(t,y_1) < x_i(t,y_2)$, defined in \eqref{E:charaeq1} for
a linearly degenerate family $i$. Let $u^\epsilon$ the wave front
solution with initial data $u(w_0^\epsilon)$, where $w_0^\epsilon$
is defined as %[E:perturb1]
\begin{equation}\label{E:perturb1}
w_0^\epsilon(x) \doteq
\begin{cases}
w(u_0(x)) & x \leq y_1 \\
w(u_0(x))+\epsilon e_i & y_1 < x \leq y_2 \\
w(u_0(x)) & x > y_2 
\end{cases}
\end{equation}
Then there exists constants $L$, $L'$, depending only on the total variation
of the initial data $u_0$,
such that for all $t \geq 0$ %[E:convsol1]
\begin{equation}\label{E:convsol1}
\int_\R \bigl| u(t,x) - u^\epsilon(t,x) \bigr| dx \leq L \epsilon
\qquad \text{and} \qquad |x_i^\epsilon(t,y_j) -
x_i(t,y_j)| \leq L' \epsilon t, \quad j = 1,2,
\end{equation}
where $x_i^\epsilon(t,y_j)$ is the position of the shock
$\sigma^\epsilon_j$ starting in $(0,y_j)$. 
\end{lemma}
\begin{proof}
The first inequality is an easy consequence of the $L^1$
continuous dependence for front tracking solutions, see \cite{BaBr}. For the
second one, note that all the shocks different from $\sigma^\epsilon$ have
size uniformly bigger than $0$, so that their position is shifted of
the order $\epsilon$. Thus the second inequality follows by standard
ODE perturbation estimates, see \cite{Br3}.
\end{proof}
An easy application of the previous lemma together with Proposition
\ref{P:finiteTV} implies that to compute $x_1(t,y_1)$ and $x_2(t,y_2)$, we
can actually consider in \eqref{E:charaeq1} a solution $\tilde u$ with
uniformly bounded total variation, so that the constant $L'$ in
\eqref{E:convsol1} is independent on the total variation of $u_0$. 

We now estimate the dependence of $x_i(t,y)$ w.r.t. $u$. %[P:shiftODE]
\begin{proposition}\label{P:shiftODE}
Let $\xi_\alpha$ be
the shift rate of the jump $\sigma_\alpha$ in $u(0,\cdot)$, and denote
with $x_i^\theta$ the solution to 
\[
\left\{ \begin{array}{l}
\dot x_i^\theta = \lambda_{i}(u^\theta(t,x_i^\theta)) \\
x_i^\theta(0) = y
\end{array} \right.
\]
where $u^\theta(t)$ is the shifted front tracking solution. Then there
exists a constant $D$ independent of the total variation of $u$ such that
%[E:shiftODE]
\begin{equation}\label{E:shiftODE}
\left| \lim_{\theta \to 0} \frac{x_i^\theta(t,y) - x_i(t,y)}{\theta} \right| \leq
D \sum_{\alpha} \Bigl| \sigma_\alpha \xi_\alpha \Bigr|.  
\end{equation}
\end{proposition}
\begin{proof}
We give a sketch of the proof, for details one can see \cite{BrGo}.

If $\epsilon$ is the size of the shock $\sigma^\epsilon$ located in
$(0,y)$, then we can apply Theorem \ref{T:shiftmap} to compute its shift
$\xi^\epsilon$: by formula \eqref{E:shifttau} we obtain %[E:approxsh1]
\begin{equation}\label{E:approxsh1}
\xi^\epsilon \sigma^\epsilon = \sum_{\alpha} \bigl( P(x+,y_\alpha) - P(x-,y_\alpha)
\bigr) \xi_\alpha.
\end{equation}
If $\theta$ is sufficiently small, then we have
\[
\xi^\epsilon = \frac{x_i^{\theta,\epsilon}(t,y)-x_i^{\epsilon}(t,y)}{\theta},
\]
where $x_i^{\theta,\epsilon}(t,y)$ is the position of the shifted shock
and $x_i^{\epsilon}(t,y)$ is its original position. Note that $(P(x+,y_\alpha)
- P(x-,y_\alpha))\xi_\alpha$ is the shift rate of the shock $\sigma^\epsilon$, after
colliding with the shocks of the Riemann problems $w_l,w_m$ and $w_m,w_r$. Their
total shift is proportional to $|\sigma_\alpha \xi_\alpha|$, and after
the interaction with $\sigma^\epsilon$, the shift of the latter is
proportional to $|\sigma^\epsilon| |\sigma_\alpha \xi_\alpha|$. 
Thus taking the limit as
$\epsilon$ tends to $0$ of \eqref{E:approxsh1}, we obtain for
$\epsilon$ sufficiently small
\[
\left| \frac{x_i^{\theta}(t,y)-x_i(t,y)}{\theta} \right| \leq D
\sum_{\alpha} \Bigl| \sigma_\alpha \xi_\alpha \Bigr|, 
\]
which implies \eqref{E:shiftODE}.
\end{proof}
We prove now the uniform Lipschitz continuity of the map $y
\longmapsto x_i(t,y)$ for all $t \geq 0$.
%[P:uniflip1]
\begin{proposition}\label{P:uniflip1}
Consider two characteristic lines $x_i(t,y_1)$ $x_i(t,y_2)$, solution
to \eqref{E:charaeq1}. There exists $C > 0$, depending only on the
system and the set $E$, such that %[E:lipdepe]
\begin{equation}\label{E:lipdepe}
{1 \over C} \leq \frac{x_i^2(t,y_2) - x_i^1(t,y_1)}{y_2 - y_1} \leq C.
\end{equation}
\end{proposition}
\begin{proof}
As in the previous proposition, let $\epsilon$ be the size of the
shock $\sigma^\epsilon(t)$ located in $(0,y)$ in Riemann coordinates.
If $\xi(t)$ is its shift rate, then for $\theta$ sufficiently small by
Theorem \ref{T:shiftmap} we obtain %[E:derivat1]
\begin{equation}\label{E:derivat1}
\frac{x_i^\epsilon(t,y+\theta \xi) - x_i^\epsilon(t,y)}{\theta}
\sigma^\epsilon(t) =
\xi^\epsilon(t) \sigma^\epsilon(t) = r_i(w_l,w_r) \bigl\langle l^i(w_l,w_r),
\sigma(0) \bigr\rangle \xi(0).
\end{equation}
In fact, by assumption, in the simplified wave patterns to compute the
shift rate of $\sigma^\epsilon$, there are no waves of the $i$-th family
different from $\sigma^\epsilon$.
Dividing by $\epsilon$ and taking the limit as $\epsilon$ tends to $0$, we obtain
\[
\frac{x_i(t,y+\theta \xi) - x_i(t,y)}{\theta} \frac{\partial}{\partial
w_i} u(t,x) = \frac{\partial}{\partial w_i}u(0,y) \bigl\langle l^i(w_l,w_r),
r_i(0,y) \bigr\rangle \xi(0),
\]
which implies 
\begin{equation}\label{E:lipdepe1}
\frac{d}{dy}x_i(t,y) = \frac{\partial u(0,y) / \partial w_i}{\partial
u(t,x) / \partial w_i} \bigl\langle l^i(w_l,w_r), r_i(0,y) \bigr\rangle.
\end{equation}
We use the fact that $\sigma^\epsilon(t) / \epsilon$ tends to
$\partial u(0,y) / \partial w_i \cdot r_i(w_l,w_r)$ as $\epsilon \to 0$. 
Since $E$ is compact,
the conclusion \eqref{E:lipdepe} follows easily.
\end{proof}
%[R:unilip1]
\begin{remark}\label{R:uniflip1}
The above proposition implies that the map $h^t_{i}$ defined in
\eqref{E:lindegmap0} is uniformly Lipschitz, independent on the total
variation of $u_0$, together with its inverse
map $(h^t_{i})^{-1}$.
\end{remark}

To end this section, we give a different proof of the following result
given in \cite{BrGo}: %[P:decay]
\begin{proposition}\label{P:decay}
If $x(t)$, $y(t)$ are the positions of two adjacent $k$-rarefaction
waves, then for some constant $\kappa > 0$ one has %[E:decay]
\begin{equation}\label{E:decay}
y(\tau) - x(\tau) \geq \kappa \tau 2^{-\nu}. 
\end{equation}
Thus for all $\tau > 0$ the total variation of the Riemann
invariant $w_{k}$ of the $k$-th genuinely nonlinear family with $N$ shocks at
$t=0$ is bounded by %[E:boundTV]
\begin{equation}\label{E:boundTV}
\TV\{w_{k}(\tau,\cdot); [a,b] \} \leq {2 (b - a) \over \kappa \tau} +
\bigl\|w_{k}\bigr\|_{L^\infty} + (N+1) 2^{1-\nu}.
\end{equation} 
\end{proposition}
\begin{figure}
\centerline{\resizebox{12cm}{6cm}{{\includegraphics{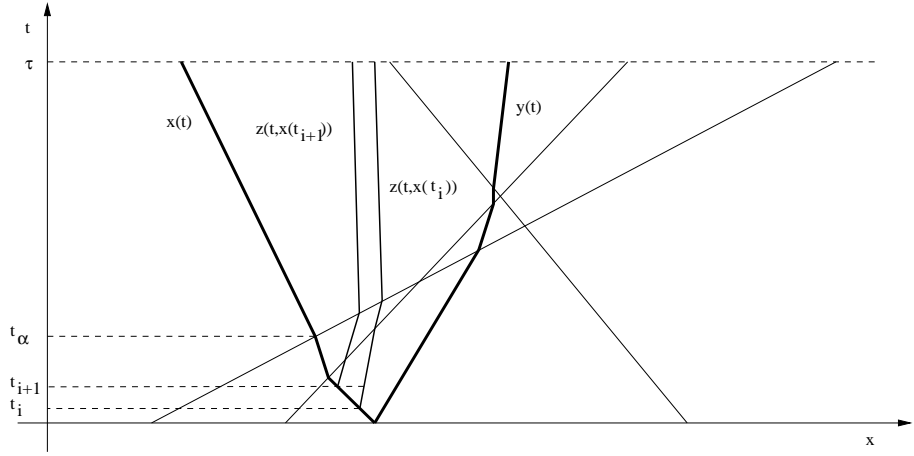}}}}
\caption{Decay of positive waves. %[Fi:decay]
}\label{Fi:decay} 
\end{figure}
\begin{proof}
Consider two adjacent $k$-rarefaction fronts $x(t)$ and $y(t)$, and
let $t_\alpha$, $\alpha=1,\dots,N$, be the interaction times of $x(t)$, $y(t)$
with other waves in the interval $[0,\tau]$. Fixed $t_i \in
(t_{\bar \alpha}, t_{\bar \alpha+1})$ for some $\bar \alpha$, let $z(t,x(t_i))$ be
the characteristic line of the $k$-th genuinely nonlinear family starting in $(t_i,
x(t_i))$ (see fig.~\ref{Fi:decay}). Assume $t_{i+1} > t_i$
sufficiently close to $t_i$ such that 
$t_{i+1} \in (t_{\bar \alpha}, t_{\bar \alpha+1})$ and $z(t,
x(t_i))$ does not collide with shocks of other
families for $t \in [t_i,t_{i+1}]$. Let $z(t,x(t_{i+1}))$ be the characteristic curve
starting in $(t_{i+1}, x(t_{i+1}))$. By the assumption
of genuinely nonlinearity, at time $t_{i+1}$ we have
\[
z(t_{i+1}, x(t_i)) - z(t_{i+1},x(t_{i+1})) \geq c (t_{i+1} - t_i)
2^{-\nu -1}, 
\]
for some constant $c$, depending only on $E$. Using Proposition
\ref{P:uniflip1}, at time $\tau$ we have %[E:chara2]
\begin{equation}\label{E:chara2}
z(\tau, x(t_i)) - z(\tau,x(t_{i+1})) \geq \frac{c}{C}
(t_{i+1} - t_i) 2^{-\nu -1}.
\end{equation}
Repeating the process, it is possible to find a countable number of
times  $t_i$ such that 
\[
\lim_{i \to -\infty} t_i = t_{\bar \alpha}, \qquad \lim_{i \to
+\infty} t_i = t_{\bar \alpha +1}, 
\]
and using \eqref{E:chara2} we get %[E:chara3]
\begin{equation}\label{E:chara3}
z(\tau,x(t_{\bar \alpha}) - z(\tau,x(t_{\bar \alpha +1})) \geq
\frac{c}{C} (t_{\bar \alpha + 1} - t_{\bar \alpha}) 2^{-\nu - 1}.
\end{equation}
Repeating the process for $y(t)$ and for all intervals $(t_{\alpha + 1},
t_{\alpha})$, we obtain \eqref{E:decay} where $\kappa =
c/C$.

The second equation follows noticing that the total
amount of positive jumps in the interval $[a,b]$ is bounded by $(1+N)
2^{-\nu} + (b-a)/\kappa \tau$.
\end{proof}

\section{Proof of the main theorem}\label{S:semicont} 
%[S:semicont]

In this section we construct the semigroup $\mathcal{S}$ on
$L^\infty(\R;E)$. In \cite{BaBr} it is shown that for all $M$, there
exists a semigroup $\mathcal{S}^M$ defined on the domain
%[E:bdTVdomain]
\begin{equation}\label{E:bdTVdomain}
D^M \doteq \Bigl\{ u : \R \mapsto E: \TV(u) \leq M \Bigr\},
\end{equation}
which is the only limit of the wave front tracking approximations
constructed in section \ref{S:setting}.
We study now the dependence of the solution on the initial data $u \in
D^M$. We consider
separately the case for genuinely nonlinear and linearly degenerate
families. %[P:gennlin1]
\begin{proposition}\label{P:gennlin1}
Consider a front tracking solution $u$, such that $u(0,\cdot)$ has $N$
jumps $\sigma_\alpha$, $\alpha=1,\dots,N$, and let $\xi_\alpha$ be
their shift rates. Given $\tau \geq 0$, denote with $\sigma_\beta$ the jumps
in the Riemann invariant $w_{k}(\tau,\cdot)$ of the $k$-th genuinely
nonlinear family. Then there exists a constant $K$,
depending only on the system and the domain $E$ such that
%[E:lipdepgnl1]
\begin{equation}\label{E:lipdepgnl1}
\sum_{\beta} \bigl| \xi_\beta \sigma_\beta \bigr| \leq K (1 + N
2^{-\nu}) \sum_{\alpha=1}^{N} \bigl| \xi_\alpha \sigma_\alpha \bigr|.
\end{equation}
\end{proposition}
\begin{proof}
The proof follows by Theorem \ref{T:shiftmap} and Proposition
\ref{P:decay}. In fact, fixed a shock $\sigma_{\bar \alpha}$, using
Theorem \ref{T:shiftmap}, we have that at time $\tau$ for a shock
$\sigma_\beta$ of the $i$-th family there exist $D'$ %[E:shiftest1]
\begin{equation}\label{E:shiftest1}
\bigl| \xi_{\beta} \sigma_{\beta} \bigr| \leq D' \bigl| \xi_{\bar
\alpha} \sigma_{\bar \alpha} \bigr|, 
\end{equation}
if the shock $\sigma_\beta$ starts on both sides of $\sigma_{\bar
\alpha}$, or, using the same estimate of Proposition \ref{P:shiftODE}, %[E:shiftest2]
\begin{equation}\label{E:shiftest2}
\bigl| \xi_{\beta} \bigr| \leq D \bigl| \sigma_{\bar \alpha} \xi_{\bar
\alpha} \bigr|, 
\end{equation}
if $\sigma_\beta$ start on one side of $\sigma_{\bar \alpha}$. Since
there is at most $1$ shocks such that \eqref{E:shiftest1} holds, and
the interval of influence is $[x_{\bar \alpha} - \hat \lambda \tau,
x_{\bar \alpha} + \hat \lambda \tau]$, using Proposition \ref{P:decay}
together with \eqref{E:shiftest1} and \eqref{E:shiftest2} we obtain
\[
\sum_\beta \bigl| \xi_\beta \sigma_\beta \bigr| \leq D' \bigl|
\xi_{\bar \alpha} \sigma_{\bar \alpha} \bigr| 
+ D \bigl| \sigma_{\bar \alpha} \xi_{\bar \alpha} \bigr| \cdot
\TV \Bigl\{ w_{k}, [x_{\bar \alpha} - \hat \lambda \tau,
x_{\bar \alpha} + \hat \lambda \tau] \Bigr\} 
\leq F ( 1 + 2^{-\nu} ) \bigl| \xi_{\bar \alpha} \sigma_{\bar
\alpha} \bigr|. 
\]
The conclusion follows the linearity of the shift differential map.
\end{proof}
Using the results of the previous section, the following result is
trivial: %[P:lindeg1]
\begin{proposition}\label{P:lindeg1}
Consider a wave front solution $u$, such that $u(0,\cdot)$ has $N$
jumps $\sigma_\alpha$, $\alpha=1,\dots,N$, and let $\xi_\alpha$ be
their shifts. Consider the equation \eqref{E:charaeq1}, whit the
eigenvalue $\lambda_i$ linearly degenerate. Fixed $\tau \geq 0$, then 
the shift $\xi_{i}$ of $x_{i}(\tau,y)$ is bounded by %[E:lipdepld1]
\begin{equation}\label{E:lipdepld1}
\bigl| \xi_{i} \bigr| \leq D \sum_{\alpha =1}^N \bigl| \xi_\alpha
\sigma_\alpha \bigr|.
\end{equation}
\end{proposition}
\begin{proof}
This is a corollary of Proposition \ref{P:shiftODE}.
\end{proof}

%We now recall that in \cite{BaBr} for every $M > 0$ there exists a
%Lipschitz continuous semigroup $\mathcal{S}^M$ of solutions of
%\eqref{E:hcl1} defined on the set $D^M$, and limit of the wave front
%tracking approximation $\mathcal{S}^\nu$. 
%For any solution $u(t,\cdot) = \mathcal{S}_t^Mu$, we recall that the map
%$h_{i_{ld}}^\tau: \R \longmapsto \R$ is defined by [E:quaslima1]
%\begin{equation}\label{E:quaslima1}
%h_{i_{ld}}^\tau(y) = x_{i_{ld}}(\tau,y),
%\end{equation}
%where $x_{i_{ld}}(t,y)$ is the solution to \eqref{E:charaeq1}:
%\[
%\left\{ \begin{array}{l}
%\dot x = \lambda_{i_{ld}}(u(t,x)) \\
%x(0) = y
%\end{array} \right.
%\]
Using the above propositions,
we can prove the following theorem: %[T:uniflipdep1]
\begin{theorem}\label{T:uniflipdep1}
Consider two initial data $u_1$ and $u_2$, and denote with
$w_{j,k}(t,\cdot)$, the $k$-th Riemann coordinate of
$\mathcal{S}^M u_j$, $j=1,2$, corresponding to the $k$-th genuinely non
linear family. Moreover, let $h_{j,i}^\tau$, $j=1,2$, the map defined in
\eqref{E:charaeq1} for the $i$-th linearly degenerate family. Then there exists a
constant $K'$, independent of $M$, such that the following
estimates hold: %[E:estima1] [E:estima2]
\begin{equation}\label{E:estima1}
\int_{\R} \bigl| w_{1,k}(t,x) - w_{2,k}(t,x) \bigr| dx
\leq K' \int_{\R} \bigl| u_1(x) - u_2(x) \bigr| dx, 
\end{equation}
\begin{equation}\label{E:estima2}
\sup_{t \geq 0, x \in \R} \bigl| h_{1,i}^t(x) -
h_{2,i}^t(x) \bigr| \leq K' \int_{\R} \bigl| u_1(x) - u_2(x)
\bigr| dx. 
\end{equation}
\end{theorem}
\begin{proof}
Consider two piecewise constant initial data $u_1^\nu$, $u_2^\nu$
in $D^{M,\nu}$, and construct a pseudo polygonal path $\gamma_0: \theta
\longmapsto u_\theta^\nu$, connecting $u_1$ and $u_2$, such that
\[
\bigl\| \gamma_0 \bigr\|_{L^1} \leq E \bigl\| u_1^\nu - u_2^\nu \bigr\|_{L^1}.
\]
We can assume that $u_\theta^\nu$ has a finite number $N$ of jumps. If
we denote with $\gamma_\tau^\nu$ the path $\theta \longmapsto
\mathcal{S}^{\nu}_\tau u_\theta^\nu$, we have by Proposition
\ref{P:gennlin1} %[E:pseudo1]
\begin{align}\label{E:pseudo1}
\bigl\| w_{2,k}^\nu(\tau) - w_{1,k}^\nu(\tau) \bigr\|_{L^1} &\leq
\Bigl\| \bigl( \gamma_\tau^\nu \bigr)_{k} \Bigr\|_{L^1} \leq K 
(1 + N 2^{-\nu}) \bigl\| \gamma_0 \bigr\|_{L^1} \\
&\leq K' (1 + N 2^{-\nu}) \bigl\| u_2 - u_1 \bigr\|_{L^1}. \notag
\end{align}
If now $\nu \to +\infty$, since $w_{j,k}^\nu(\tau)$ converges to
$w_{j,k}(\tau)$, we obtain \eqref{E:estima1}. Since this
estimate does not depend on the number of initial jumps $N$, we can
extend it uniformly on $D^M$. 

Using the same pseudo polygonal path, in a similar way we can prove that 
\[
\bigl| x_{2,i}^\nu(\tau,y) - x_{1,i}^\nu(\tau,y) \bigr| \leq K'
\bigl\| u_2 - u_1 \bigr\|_{L^1}.
\]
This shows that $x_{i}^\nu(\tau,\cdot)$ converges uniformly to the
solution $x_{i}(\tau,\cdot)$ as $\nu \to +\infty$ and $u^\nu \to
u$. It also implies that
\[
\bigl| x_{2,i}(\tau,y) - x_{1,i}(\tau,y) \bigr| \leq K'
\bigl\| u_2 - u_1 \bigr\|_{L^1},
\]
This concludes the proof.
\end{proof}

We can now define $\mathcal{S}$ on the domain $L^\infty(\R;E)$:
%[D:solution1]
\begin{definition}\label{D:solution1}
For all $u \in L^\infty(\R,E)$, let $u^M \in D^M$ be such that
%[E:nearu1]
\begin{equation}\label{E:nearu1}
\lim_{M \to +\infty} u^M = u \quad \text{in }L^1_{\text{loc}}.
\end{equation}
Define $\mathcal{S}_t u$ as %[E:semigroT1]
\begin{equation}\label{E:semigroT1}
\mathcal{S}_t u = \lim_{M \to +\infty} \mathcal{S}_t^{M} u,
\end{equation}
where the limit is in $L^1_{\text{loc}}$.
\end{definition}
It is easy to prove that the right hand side of \eqref{E:semigroT1} is
a Cauchy sequence in every compact set $[a,b]$: 
in fact, using the finite speed of propagation, we can consider $u$
with compact support $[a - \hat \lambda t, b + \hat \lambda t]$.
For the components $w_{k}$ of the $k$-th genuinely nonlinear family, it follows
directly from \eqref{E:estima1}, while for a linearly degenerate
component $w_{i}$, let $\tilde w$ be a Lipschitz
continuous function such that
\[
\int_{\R}  \bigl| w_{i}(0,x) - \tilde w(x) \bigr| dx \leq \epsilon.
\]
By Theorem \ref{T:uniflipdep1} we have for $u_1,u_2 \in D^M$ such
that $\|u - u_i\|_{L^1} < \delta$, $i=1,2$,
\[
\sup_{t \geq 0, x \in \R} \bigl| h_{1,i}^t(x) -
h_{2,i}^t(x) \bigr| < K' \delta,
\]
and it follows by easy computations that
%[E:nonunif1]
\begin{align}\label{E:nonunif1}
\bigl\| w_{1,i}(t) - w_{2,i}(t) \bigr\|_{L^1}
&\leq \bigl\| w_{1,i}(t) - \tilde w \circ
\bigl(h^t_{1,i}\bigr)^{-1} \bigr\|_{L^1} + \bigl\|
w_{2,i}(t) - \tilde w \circ \bigl( 
h^t_{2,i} \bigr)^{-1} \bigr\|_{L^1} + \\
& \qquad \qquad \bigl\| \tilde w \circ
\bigl(h_{1,i}^{t} \bigr)^{-1} - \tilde w \circ \bigl( 
h^t_{2,i} \bigr)^{-1} \bigr\|_{L^1} \notag \\
&\leq C \bigl\| w_{1,i}(0) - \tilde w \bigr\|_{L^1} +
C \bigl\| w_{2,i}(0) - \tilde w \bigr\|_{L^1} 
+ L (b - a) G \bigl\|u_2 - u_1\bigr\|_{L^1} \notag \\
&\leq 2 C (\epsilon + \delta) + L (b - a) G \delta, \notag
\end{align}
where $L$ is the Lipschitz constant of $\tilde w$. %and 
%\[
%\int_{a - \hat \lambda t}^{b + \hat \lambda t} \bigl| u_1(x) - u_2(x)
%\bigr| dx \leq \delta.  
%\]
This shows that $w_{i}^M(t)$ is a Cauchy sequence for all $t \geq
0$, because the right hand side of \eqref{E:nonunif1} can be made
arbitrarily small. We can now prove the main theorem: %[T:final1]
\begin{theorem}\label{T:final1}
The semigroup $\mathcal{S} : [0,+\infty) \otimes L^\infty(\R;E)
\longmapsto L^\infty(\R;E)$ defined in \eqref{E:semigroT1} is the 
only continuous semigroup on $L^\infty(\R;E)$ such that the following
properties are satisfied:
\begin{itemize}
\item[i)] for all $\bar u^n, \bar u \in L^\infty(\R;E)$, $t_n,t \in [0,+\infty)$,
with $\bar u_n \to \bar u$ in $L^1_{\text{loc}}$, $|t - t_n| \to 0$ as
$n \to +\infty$, %[E:contin1] 
\begin{equation}\label{E:contin1}
\lim_{n \to +\infty} \mathcal{S}_{t_n} \bar u_n = \mathcal{S}_{t} \bar
u \quad \text{in} \ L^1_{\text{loc}}; 
\end{equation}
%\item[i)] for all $u^n,u \in L^\infty(\R;E)$, $t_n,t \in [0,+\infty)$,
%with $\|u - u_n\|_{L^1_{\text{loc}}} \to 0$, $|t - t_n| \to 0$ as $n \to +\infty$, %[E:contin1]
%\begin{equation}\label{E:contin1}
%\lim_{n \to +\infty} \bigl\| \mathcal{S}_t u - \mathcal{S}_{t_n} u_n
%\bigr\|_{L^1_{\text{loc}}} = 0; 
%\end{equation}
\item[ii)] each trajectory $t \mapsto \mathcal{S}_t u_0$ is a weak entropic
solution to the Cauchy problem %[E:hcl2]
\begin{equation}\label{E:hcl2}
\left\{ \begin{array}{l}
u_t + f(u)_x = 0 \\
u(0,x) = u_0(x)
\end{array} \right.
\end{equation}
with $u_0 \in L^{\infty}(\R;E)$;
\item[iii)] if $u_0$ is piecewise constant, then, for $t$ sufficiently
small, $\mathcal{S}_t u_0$ coincides with the function obtained by
piecing together the solutions of the corresponding Riemann problems.
\end{itemize}
\end{theorem}
\begin{proof}
The statement follows easily, since we proved that $\mathcal{S}_t u$
is the unique limit of wave front approximations, and for data with
bounded total variation we can apply the results in \cite{BaBr}. 
\end{proof}
%[R:finalrem]
\begin{remark}\label{R:finalrem}
Note that what we also proved that the characteristic equation
\eqref{E:charaeq1} is well posed for $L^\infty$ data: the
solution $x_{i}(t,y)$ is Lipschitz continuous w.r.t. both variables. 
This is not trivial, since even for $2 \times 2$ systems not in
conservation form the dependence is H\"older continuous, while for
general $n \times n$ the solution does not exists. 

Note moreover that semigroup $\mathcal{S}$ is continuous, but not uniformly
continuous. However if the initial data takes values is a compact set
of $L^1 \cap L^\infty$, then the semigroup becomes uniformly
continuous. This extend the Lipschitz continuity when the initial data
have bounded total variation.
%: in fact for the genuinely nonlinear components one applies the decay
%of the total variation, while for the others the compactness follows
%from the Riesz-Fr\'echet-Kolmogorov theorem on compact sets of
%$L^p$ (see \cite{Brez}).
\end{remark}

\end{document}